\documentclass[sn-mathphys]{sn-jnl}

\usepackage[utf8]{inputenc}
\usepackage{graphicx}%
\usepackage{multirow}%
\usepackage{amsmath,amssymb,amsfonts}%
\usepackage{amsthm}%
\usepackage{mathrsfs}%
\usepackage[title]{appendix}%
\usepackage{xcolor}%
\usepackage{textcomp}%
\usepackage{manyfoot}%
\usepackage{booktabs}%
\usepackage{algorithm}%
\usepackage{algorithmicx}%
\usepackage[capitalise,noabbrev]{cleveref}%
\usepackage{algpseudocode}%
\usepackage{listings}%
\usepackage{soul}

\usepackage{tikz}
\usetikzlibrary{arrows,petri,topaths,shapes,decorations,fit,positioning}
\usepackage{pifont}
\usetikzlibrary{patterns}
\usepackage{pgfplots}
\usepackage{multicol}
\usepackage{subcaption}
\usepackage{changepage}

\usetikzlibrary{patterns.meta}

\input{command.sty}

\begin{document}

\title{Fairness by design in shared-energy allocation problems}



\author*[1,2]{\fnm{Zo\'e} \sur{Fornier}}\email{zoe.fornier@enpc.fr}

\author[2]{\fnm{Vincent} \sur{Lecl\`ere}}\email{vincent.leclere@enpc.fr}

\author[3,4,5]{\fnm{Pierre} \sur{Pinson}}\email{p.pinson@imperial.ac.uk}

\affil*[1]{\orgname{METRON},\country{France}}

\affil[2]{\orgname{CERMICS}, \country{France}}

\affil[3]{\orgname{Imperial College London}, \country{United Kingdom}}
\affil[4]{\orgname{Technical University of Denmark}, \country{Denmark}}
\affil[5]{\orgname{Halfspace}, \country{Denmark}}


\abstract{
This paper studies how to aggregate prosumers (or large consumers) and their collective decisions in electricity markets, with a focus on fairness. 
Fairness is essential for prosumers to participate in aggregation schemes. Some prosumers may not be able to access the \change{electricty} market directly, even though it would be beneficial for them. Therefore, new companies offer to aggregate them and promise to treat them fairly. This leads to a \emph{fair resource allocation} problem.

We propose to use \emph{acceptability constraints} to guarantee that each prosumer gains from the aggregation. 
Moreover, we aim to distribute the costs and benefits fairly, taking into account the multi-period and uncertain nature of the problem. 
Rather than using financial mechanisms to adjust for fairness issues, we focus on various objectives and constraints, within decision problems, that achieve \emph{fairness by design}. We start from a simple single-period and deterministic model, and then generalize it to a dynamic and stochastic setting using, \eg stochastic dominance constraints.
    
   
}




\maketitle

\section{Introduction}\label{sec:introduction}

Many domains, such as telecommunication networks, healthcare, disaster management, and energy-sharing systems, require fairness as a key criterion. 
However, fairness is not easy to define or implement, as it can have different meanings and implications in different contexts. 
\change{Nevertheless}, mathematical models that address real-world problems should not ignore fairness, even if it adds complexity to the problem. 
In this paper, we investigate various methods to incorporate fairness in a multi-agent problem. Specifically, we apply fairness to the problem of aggregating prosumers, who are both \change{electricity} producers and consumers, in the \change{electricity} market.

We focus on electric energy management application, where the aggregation of prosumers is becoming more relevant due to the increasing number of prosumers. Renewable energy generation capacities are becoming more affordable and effective, as renewable energy investments are rising (19\% in 2022, according to a report by \cite{IRENA} on global trends in renewable energy). This enables smaller prosumers, such as medium-sized industries, to invest in onsite energy generation and storage. However, prosumers are usually too small to access the electricity market directly, so some companies offer to aggregate them in \change{electricity} markets.
\new{For example, CPower  is an American company that aggregates a total of $2.000$ MW of power \cite{CPOWER}. We refer to \cite{CARREIRO20171160} for an extensive review on aggregators and their role in electricity markets.}

Those aggregators can be external entities responsible for every prosumer energy transfers.
In this case, there is a necessity to think of how the aggregation affects the participants to ensure a fair allocation of benefits.
This is highlighted in a report \citep{eurelectric} on designing fair and equitable market rules for demand response aggregation, 
published by the association representing the common interests of the European electricity industry, Euralectric.
Indeed, there is a practical need to guarantee that each prosumer benefits from staying in the aggregation.
Further, 
prosumers need to feel like they are not being disfavored compared to others,
leading the aggregator to choose a solution with a fair allocation of benefits.

In the literature, one distinguishes two main approaches in handling fairness: 
solve the problem efficiently and then reallocate the benefits \citep{YHS2021,WZW+2019a,YFLM2023}; 
or change the objective function in order to get a fair solution \citep{Hoo}. 
In the first approach, we model a multi-agent problem with a utilitarian objective, i.e., we optimize the aggregated objectives of agents.
Then, a protocol is implemented to reallocate the benefits among agents.
For example, \textit{Shapley values} \citep{P-295} assess the marginal contribution of each agent in the group and determine their fair share.
The second approach prioritizes fair solutions through the modeling by changing the objective function.
\change{
We refer to \cite{Hoo} for a comprehensive overview and guidelines on selecting an appropriate objective function to reflect fairness.}
The two most studied objective functions are \textit{the minimax objective} \citep{Raw71}, which optimizes the least well-off agent's objective, and \textit{the proportional objective} \citep{Nash}, derived from Nash's \textit{bargaining solution}, which optimizes the logarithmic sum of agents' objectives.
\change{Note that this approach, through the objective function, means that the well-being of different agents are compared through a single value.}

However, these approaches present some limitations.
On the one side, the proportional and minimax approaches focus merely on the objective function and not decisions.
This can be a problem, as in some applications there can be different characteristics which are valuable. 
For example in an energy contract, both the flexibility and the volume of energy traded are important features.
Thus, it is hard to take into account both of them when the quality of a solution is determined by a single value.
On the other hand, 
post-allocation distributions of benefits are not adapted to problems that are formulated over long periods of time, such as contracts in \change{electricity} markets. 
Indeed, those approaches require to solve the whole problem before allocating costs. 
Then, it is impractical in most cases to expect each agent to wait until the problem's completion, which could span several months or years, 
to receive their fair share.
Furthermore, given the inherent uncertainties linked to most problems, we also want our approach to hold in a stochastic framework.
Then, fairness criteria must be redefined considering utility distributions and associated risks over time.

In this paper, we introduce various strategies for integrating fairness considerations into optimization problems. 
Our primary focus is what we refer to as \textit{fairness-by-design}.
Instead of relying on \textit{ex post} redistribution, \change{as is usual in} game theory \citep{P-295}, we can establish a degree of fairness directly within the model. 
\new{Our main contribution is to provide a framework and tools to accommodate fairness into mathematical models, in particular in the context of prosumer aggregation. 
What sets our approach apart is the extension of this framework to dynamic and stochastic settings, allowing for risk-averse and time-consistent guarantees.}

More specifically, we present two key elements for achieving fair allocation in an aggregation. 
\change{First, we model fair cost allocation through an operator ordering the costs of the different prosumers. 
For the choice of this operator, we present three traditional approaches (\textit{utilitarian, proportional} and \textit{minimax}).
Additionally, we propose \textit{acceptability constraints}
that ensure each agent's outcome improves in a predefined sense within the aggregation. 
In their simplest form, these acceptability constraints correspond to individual (or self) rationality in game theory, ensuring each agent benefits from being in the group.
We then extend the problem to a dynamic framework where decisions are made sequentially over time. 
In this context, agent's cost are multidimensional, the acceptability (or individual rationality) constraints thus need to choose a (partial) order. We discuss a few relevant partial order choices.
Similarly, in a stochastic framework, agents cost are random variables, and we discuss relevant stochastic orders.}
Compared to \cite{GKW2023}, who propose a risk-averse stochastic bargaining game, 
our approach handles uncertainties through the objective function but also acceptability constraints. 
This enables us to consider various aspects of the impact of uncertainties on the problem.
As a result, our proposed model is well-suited for addressing inherent uncertainties within multistage stochastic programs, enhancing its practical applicability. 
Finally, we assess these different strategies on a toy model \new{ where we aggregate $4$ electricity consumers to access the day-ahead market.} 
We discuss each modeling choice consequences.

The remainder of the paper is organized as follows. 
In \cref{sec:fairness}, we delve into definitions of fairness and its integration into optimization models.
We propose, in \cref{sec:problem}, to model prosumers aggregation with acceptability constraints and a fair objective function.
\new{We then illustrate the introduced framework on a toy model in \cref{subsec:deterministic_illustration}}.
\Cref{sec:dynamic} expands the notion of acceptability into the dynamic framework,
while \cref{sec:stochastic} adapts acceptability and fairness to the stochastic framework.

\new{
\subsubsection*{Notations}

To facilitate understanding, we go through some notations used in this paper. 
We denote $[N] := \na{1, \dots, N}$ the set of non-null integers smaller than T . 
Accordingly, $X_{[n]}$ denote the collection $\na{X_i}_{i\in[n]}$. 
Random variables are denoted in bold characters, with their realization is normal font.
The $\sigma -$algebra generated by $\na{\xi_\tau}_{\tau \in [t]}$ is denoted $\sigma(\xi_{[t]})$.
Finally, in this paper, the term operator always refers to a mathematical operator.
}

\section{\change{Fairness in the literature}}
\label{sec:fairness}
In the Oxford Dictionary, fairness is defined as \textit{the quality of treating people equally or in a way that is reasonable}.
The definition is simple but subjective.
Is treating people equally, regardless of any token of individuality, considered fair in society?
Furthermore, what does it mean to be reasonable? 
Whatever take we have on fairness is necessarily subjective and context-dependent (see \cite{Kon2003} for a philosophical analysis of fairness).
In this section, we give a general overview of how fairness is defined and modeled across the scientific literature, while linking it to our energy application. 
Bear in mind that each approach on fairness adopts a specific definition of fairness, which is not consensual.

\subsection{Modeling and accommodating fairness}
\label{subsec:objectives_literature}

One of the main challenges facing fairness is the allocation of resources between agents.
This naturally falls into the scope of Game Theory,
where each of the $N$ individuals (or players) is modeled with a utility function whose actual value depends on the actions of all players.
For a given set of actions, we obtain a utility vector, denoted $u:=(u_1, \dots, u_N)$, representing the utility of every agent $u_i$. 
A utility vector is then said to be \emph{fair} if it satisfies a set of properties that might vary from one specific fairness definition to another.
\new{Among them, \emph{individual rationality}, which ensures that every individual is better off in the aggregation, and therefore accepts to be a part of it, is often required.}

\temporary{Game Theory: Nash solution}
In a seminal contribution \citep{Nash}, John Nash introduced the bargaining problem
where two \new{rational} agents, allowed to bargain, try to maximize the sum of their utilities. 
\new{If agents are rational, individual rationality must be ensured. 
This is modeled using a \textit{disagreement point} which represents the \change{outcome} obtained by players if they cannot reach an agreement.}
For agents to cooperate, they must agree on properties a utility vector, $u_{[N]}$, should satisfy to be admissible. 
Nash proposed four axioms to constitute this agreement:
\textit{Pareto optimality} -- we cannot improve the utility of one agent without decreasing another's utility;
\textit{Symmetry} -- applying the same permutation to two utility vectors does not change their order;
\textit{Independence of irrelevant alternatives} -- if a utility vector is the optimal utility vector within the feasible set, it remains so if the set is reduced;
\textit{Scale invariance} -- applying affine transformations to the utility vector does not change the social ranking. 
Nash showed that, under some assumptions including convexity and compactness of feasible utility vectors,
there exists a unique utility vector satisfying those axioms. 
This unique utility vector is regarded in the literature as a viable option when seeking fairness.
It has been demonstrated that under convexity of the feasible set, it can be obtained by maximizing the product of utilities (\cite{Nas1953,Mut1999}), and thus by maximizing a logarithmic sum of utilities:
\begin{align*}
    \max_{u \in \cU} & \quad \sum_{i=1}^N \log(u_i - d_i), 
\end{align*}
where $u\in \cU$ is a feasible utility vectors among $N$ players, and $d$ the disagreement point.
This approach is referred to as \textit{proportional fairness}.
Some papers criticized the \textit{Independence of irrelevant alternatives} for having undesirable side effects. 
To overcome those issues, \cite{Kalai} proposed to replace it with a \textit{monotonicity} axiom, resulting in another unique utility vector
and a slightly different vision on fairness.

\temporary{Game Theory: Cooperative games and shapley values}
In contrast to bargaining games, \textit{ cooperative games} study games in which coalition formation is allowed:
see \cite{Osborne1994} for a complete introduction.
In this theory, it is assumed that players can achieve superior outcomes by cooperating \delete{rather than working against each other.}
Players must establish their common interest and then work together to achieve it, 
which requires information exchanges.
\new{In \textit{transferable utility games}, payoffs are given to the group which then divides among players through a post-allocation scheme.
}
In \cite{P-295}, Shapley studied a class of functions that evaluate the participation of players in a coalition. 
Considering a set of axioms (symmetry, efficiency and law of aggregation), Shapley showed that there exists a unique value function satisfying those axioms.
He derived an explicit formula to compute the value of a player $i$ in a cooperative game with a set $N$ of players: 
\begin{align*}
    \phi_i(v) = \sum_{S \subset N \setminus \na{i}} \binom{|N|-1}{|S|}^{-1} \, \bp{v(S \cup \na{i}) - v(S)},
\end{align*}
where $v(S)$ gives the total expected sum of payoffs the coalition $S$ can obtain.
The values obtained $\na{\phi_i(v)}_{i \in N}$ are called Shapley values.
They are considered a fair redistribution of gains in the group. 
However, they are very hard to compute in practice (as the size of the problem grows, those values are not computable).

\temporary{The Minimax approach}
A different approach 
was introduced by John Rawl in \cite{Raw71}:
assuming that a group of individuals has no idea of their rank or situation in society,
they will agree on a social contract aiming at maximizing the well-being of the least well-off. 
If the agents possess distinct characteristics, 
it might be difficult to compare them and ensure equitable treatment among them.
This approach to fairness is often referred to as \textit{minimax fairness}, as this amounts to optimizing the worst objective among agents.

\temporary{The price of fairness}
When fairness is considered in the problem (through the objective or constraints), it comes at a price: a fair solution might not be the most efficient one. 
Indeed, many articles try to find a balance\change{, or trade-off,} between efficiency (have the best objective possible) and fairness (have a fair solution). 
In \cite{BFT2011}, the authors established bounds on the price of fairness for two approaches, \change{previously introduced- proportional fairness and minimax fairness-} in resource allocation problems among self-interested players.

\temporary{Reference for a more complete overview}
In this section, we referred to work that lay the foundations of fairness modeling in mathematics. 
In the following section, we present some applications of aggregations and the way fairness is considered or evaluated. 

\subsection{Applications of fairness in the literature}

\temporary{For simplicity, the utilitarian approach is usually chosen.}
In this paper, we focus on a \textit{by-design} approach, meaning that fairness is already accommodated in an optimization model. 
Although fairness is commonly recognized as crucial, in most articles the approach adopted derives from act utilitarianism: 
one should at every moment promote the greatest aggregate happiness, which consists in maximizing social welfare regardless of individual costs.
For example, in \cite{XIAO2020115159}, the authors studied an aggregator in charge of multiple agents within a power system.
They optimized the total revenue of the aggregation without considering the impact on each agent individually. 
In \cite{MP2019}, a prosumers' aggregator can focus on different indicators (import/export costs, exchange with the system operator, peak-shaving services etc.) to optimize its trades with the energy market, and the trades between prosumers.
The indicator to focus on must be agreed on by the prosumers.
The authors gave a sensitivity analysis of the parameters of the problem to determine what would increase the social acceptability of such an aggregation system. 
However, the model is utilitarian as it does not consider the allocation of costs among agents.

\temporary{Coalitional games.}
Other papers have proposed to first optimize the problem and then \change{handle fairness through benefit post-allocation schemes.}
\new{One way to deal with post-allocation is to }
model the aggregation as a coalitional game.
This is the case of \cite{FSLB2015}, where the authors studied a risk-averse renewable-energy multi-portfolio problem. 
In order to get a fair and stable allocation of profits, they chose \textit{the Nucleolus} approach which finds a vector utility that minimizes the incentive to leave the aggregation for the worst coalition.
In particular, this solution is in \textit{the core} of the game, meaning every players gains from staying in the grand coalition.
Similarly, in \cite{YHS2021}, the authors studied a group of buildings with solar generation that mutually invest in an ESS. 
The approach is to, first, optimize the problem formulated as a two-stage stochastic coalition game.
Then, a fair reallocation of costs is determined by computing the nucleolus allocation which minimizes the minimal dissatisfaction of agents.

\temporary{Post allocation schemes}
Some papers \change{propose different methods to elaborate} post-allocation schemes.
For example in \cite{YFLM2023}, the authors studied the joint participation of wind farms with a shared energy storage. 
The solution is found by first solving a two-stage stochastic program, 
and then reallocating the lease cost among users in a proportional scheme.
They chose to make a wind farm pay depending on its increase of revenue after using the energy storage leasing service.
In \cite{WZW+2019a}, the authors valued cooperation in their model, which is another way to look at cost redistribution.  
They considered an aggregator which participates in capacity and energy market for a number of energy users.
In their model, the aggregator is not in charge of the users decisions but of the trades with the energy market, 
therefore he must incentivize users to deviate from their optimal scheduling for minimizing total revenue. 
They proposed to solve an asymmetric Nash bargaining problem to determine the incentizing costs.
In another approach, \change{
the authors solved a multi-portfolio problem with fairness considerations in \cite{IT2014}. 
Instead of splitting the market impact costs in a pro-rata fashion, they introduced charging variables, constrained to satisfy some properties, that are optimized in the model.
This approach amounts to having transfer variables, which we avoid in this paper, as they may raise privacy and trust concerns in practical application. 
Instead, we simplify the approach by designating the aggregator as the sole entity with complete information on the problem, \change{which pays directly agents depending on their actions}.
}

\temporary{Fairness in the constraints}
Typically, fairness is dealt with through the objective function, or in a post-allocation scheme. 
However, some researchers proposed constraints to ensure fairness.
For example in \cite{AKY2022}, the authors constrained the allocation feasibility set for a resource allocation problem. 
They introduced a welfare function dominance constraint: 
the admissible set of social welfare functions must dominate a referenced one.
Then, with a utilitarian objective, a trade-off between fairness and efficiency is obtained.
An alternative approach, proposed in \cite{Oh2022}, is to bound a fairness indicator. 
The authors studied the energy planning of multiple agents over a virtual energy storage system (VESS), where energy dispatch is managed by an aggregator. 
They introduced two fairness indicators depending on the energy allocation, and added constraints bounding them in a utilitarian model. 
Then, they compared the results with a minimax approach, where they optimize the minimal fairness indicator over agents. 

\temporary{Stochastic considerations}
In many cases, uncertainties are inherent to the problem. 
If multiple articles have dealt with uncertainties, they rarely have a stochastic take on fairness. 
For example, in both \cite{YFLM2023} and \cite{YHS2021}, the authors solved their problem with a two-stage program and then redistributed the costs fairly after uncertainty realization. 
Thus, there is no stochastic policy for fair redistribution. 
Other articles accommodated risk-averse profiles to game theory approaches. 
In \cite{GKW2023}, the authors studied a risk-averse extension of the Bargaining Problem. 
They adapted Nash bargaining axioms to constrain the feasible utility vectors depending on the risk profile of players.

\section{A shared-resource allocation problem in the context of a prosumer aggregator}
\label{sec:problem}

\change{We present here a general framework where a so-called \emph{aggregator}
aggregates independent agents' needs (industrial prosumers, residential units, virtual power plants\dots) 
and makes economic transactions for the collective. }
To make aggregation contracts attractive to \change{agents},  we encounter two distinct challenges:
first, each \change{agent} needs to find the contract \emph{acceptable}, ensuring that each agent derives substantial benefits from the aggregation;
second, the decisions made by the aggregator, leading to benefits or losses for each agent, should be made fairly. 
Recall that, for practical reasons, we do not allow money transfers between agents.
\new{Finally, to align with standard optimization frameworks, we aim to minimize the costs of agents and thus consider them as buyers.}

In the following, \cref{ssec:market_structure} formalize the setting, 
\cref{subsec:fairness_deterministic} explore various objective functions that model fair decisions, and finally
\cref{subsec:acceptibility_deterministic} introduce acceptability constraints.

\subsection{Prosumers and market structure}
\label{ssec:market_structure}


We denote by $\decisionProsumer{\indexProsumer} \in \constraintProsumer{i}$ the set of state and decision variables modeling \change{an agent} $\indexProsumer$.
The technical constraints \change{proper to agent} $\indexProsumer$ are represented through feasible set $\constraintProsumer{i}$, 
while \change{external constraints (for instance market exchanges)},
common to all agents, are represented with feasible set $\constraintMarket$.
Finally, each prosumer wants to \change{minimize} a cost function $\objectiveProsumer{\indexProsumer}:\constraintProsumer{\indexProsumer} \rightarrow \RR$, yielding the model $\modelProsumer{\indexProsumer}$.
\change{Note that $\modelProsumer{\indexProsumer}$ can model problems in various contexts. }
\new{In \cref{subsec:deterministic_illustration}, we present the particular application of this framework to prosumers aggregation on energy markets.}

We now consider an aggregator in charge of $\setProsumer$ agents,
we denote $\decisionProsumer{ } := (\decisionProsumer{\indexProsumer})_{\indexProsumer \in [\setProsumer]}$.
The aggregator in problem $(A)$, 
\change{aggregates agents' decisions into $h(\decisionProsumer{1},\dots, \decisionProsumer{\setProsumer})$ to satisfy external constraints $\cM$ (see \eqref{eq:binding_constraint}).
Further, the physical constraint of each agent must be conserved (see \eqref{eq:aggregator-constraint-physic}),
while the external constraints bind all agents' decisions.}
Finally, on the one hand, constraint \eqref{eq:abstract_acceptability_constraint} ensures that the cost of an agent $\indexProsumer$ is 
within an acceptable set $\acceptableSet{\coefficientAcceptability}{\indexProsumer}$ they have agreed on prior to optimization.
On the other hand, $\operatorAggregation{\setProsumer}$ is the \operatorName that computes the objective of the aggregator considering the $\setProsumer$ objective functions of all agents. 
Depending of the choices of the acceptability sets $\acceptableSet{\coefficientAcceptability}{\indexProsumer}$ and the \operatorName $\operatorAggregation{\setProsumer}$, discussed, respectively, in \cref{subsec:fairness_deterministic} and \cref{subsec:acceptibility_deterministic}, 
We have obtained different approaches to the shared resource allocation problem.

\begin{subequations}
\begin{alignat}{5}
    \modelProsumer{\indexProsumer} \;  \qquad \Min_{\decisionProsumer{\indexProsumer}} \quad & \objectiveProsumer{\indexProsumer}(\decisionProsumer{\indexProsumer})  & 
    \qquad
    \modelAggregator \; \quad \Min_x \quad & \operatorAggregation{\setProsumer}((\objectiveProsumer{\indexProsumer}(\decisionProsumer{\indexProsumer}))_{\indexProsumer \in [\setProsumer]}) 
    \\
    \text{s.t.} \quad & \decisionProsumer{\indexProsumer} \in \constraintProsumer{\indexProsumer} &  
    \text{s.t.} \quad & \decisionProsumer{\indexProsumer} \in \constraintProsumer{\indexProsumer} & \quad \forall \indexProsumer \in [\setProsumer] \label{eq:aggregator-constraint-physic} \\
    & \decisionProsumer{\indexProsumer} \in \constraintMarket. & & \change{h(\decisionProsumer{1},\dots, \decisionProsumer{\setProsumer})} \in \constraintMarket &\label{eq:binding_constraint}  \\
    & & & \objectiveProsumer{\indexProsumer}(\decisionProsumer{\indexProsumer}) \in \acceptableSet{\coefficientAcceptability}{\indexProsumer} & \forall \indexProsumer \in [\setProsumer]. \label{eq:abstract_acceptability_constraint}
\end{alignat}
\end{subequations}

\new{We assume that the aggregation can decide that it is optimal for agents to operate independently \ie if $\decisionProsumer{\indexProsumer}$ is an optimal solution of $\modelProsumer{\indexProsumer}$, then ($\decisionProsumer{1}$, \dots, $\decisionProsumer{\setProsumer}$) is an admissible solution of $\modelAggregator$.}

\subsection{Fair cost aggregation}
\label{subsec:fairness_deterministic}

Assuming that all agents have agreed to participate in the aggregation (we discuss acceptability in \cref{subsec:acceptibility_deterministic}),
we focus on the way the aggregator operates to allocate
aggregation benefits among prosumers.

The most natural and efficient method is the so-called \textit{utilitarian approach}:
\begin{subequations}
\begin{align}
    \objectiveUtilitarian{\setProsumer}((\objectiveProsumer{\indexProsumer}(\decisionProsumer{\indexProsumer}))_{\indexProsumer \in [\setProsumer]}) = \sum_{\indexProsumer \in [\setProsumer]} \objectiveProsumer{\indexProsumer}(\decisionProsumer{\indexProsumer}). \label{obj:utilitarian}
\end{align}
This approach aims to minimize total costs independently from the distribution of costs among prosumers: fairness is set aside. 
Indeed, in case of heterogeneity of the objective functions, 
\new{
it is possible that one of the objective function $\objectiveProsumer{\indexProsumer}$ dominates the others, \ie $$\objectiveProsumer{\indexProsumer}(\decisionProsumer{\indexProsumer}) \geq \objectiveProsumer{k}(\decisionProsumer{k}), \qquad \forall \decisionProsumer{\indexProsumer} \in \constraintProsumer{\indexProsumer}, \quad \forall \decisionProsumer{k} \in \constraintProsumer{k},$$
in which case
}
all efforts of the aggregation are focused on minimizing the dominant objective function. 
A possibility that falls out of the scope of this paper (see \cref{subsec:objectives_literature}) is to solve $\modelAggregator$ and then reallocate resources with a fair scheme or put money transfers in place.
We study alternative \operatorNames that ensure fair allocation for various fairness definitions.

First, we consider the \emph{proportional approach} based on Nash bargaining solutions (see \cref{subsec:objectives_literature}).
For this approach, we consider the set of reachable (dis)utilities
$\cL = \ba{(\objectiveProsumer{1}(\decisionProsumer{1}), \dots, \objectiveProsumer{\setProsumer}(\decisionProsumer{\setProsumer}) \;|\; 
\decisionProsumer{\indexProsumer} \in \constraintProsumer{\indexProsumer}, \, \forall \indexProsumer \in [\setProsumer], \;  \matrixMarketProsumer{\indexProsumer} \decisionProsumer{\indexProsumer} \in \constraintMarket
}$, and set the optimal values of $\modelProsumer{\indexProsumer}$, $\valueProsumer{\indexProsumer}$ \change{\footnote{We implicitly assume here that either there is a unique solution, or that we have defined a way to select a solution among the set of optimal solutions.}}, 
as the chosen disagreement point.
Then, \citet{Nash} introduces a set of axioms that must respect a fair distribution of (dis)utilities, 
and show that, if 
$\cL$ is convex and compact, there exists a unique (dis)utility vector satisfying those axioms.
Furthermore, it is proven that Nash's distribution is obtained 
by maximizing the sum of logarithmic utilities.
For our problem, it corresponds to using the \operatorName:
\begin{align}
    \objectiveProportional{\setProsumer}((\objectiveProsumer{\indexProsumer}(\decisionProsumer{\indexProsumer}))_{\indexProsumer \in [\setProsumer]}) := \,  - \sum_{\indexProsumer \in [\setProsumer]} \log( v^\indexProsumer - \objectiveProsumer{\indexProsumer}(\decisionProsumer{\indexProsumer})). \label{obj:nash}
\end{align}
Note that this approach tends to act in favor of smaller participants.
Indeed, increasing a small cost improvement is preferred to increasing an already large cost improvement.

Finally, Rawls' theory of justice leads to the \textit{minimax approach} favoring the least well-off. 
Here, the operator we obtain is:
\begin{align}
    \objectiveMinimax{\setProsumer}((\objectiveProsumer{\indexProsumer}(\decisionProsumer{\indexProsumer}))_{\indexProsumer \in [\setProsumer]}) := \max_{\indexProsumer \in [\setProsumer]} \; \objectiveProsumer{\indexProsumer}(\decisionProsumer{\indexProsumer}).  \label{obj:minimax}
\end{align}
\change{For similar reasons to the utilitarian approach, this method may not be adequate for heterogeneous agents as it only focuses on minimizing the dominant objective function. }
To address this issue, we quantify an agent's well-being by looking at the proportional savings he makes in the aggregation.
Then, applying Rawls' principle\change{, we minimize the maximum proportional costs over agents, and} we obtain the following agent operator:

\begin{align}
    \objectiveMinimaxSavings{\setProsumer}((\objectiveProsumer{\indexProsumer}(\decisionProsumer{\indexProsumer}))_{\indexProsumer \in [\setProsumer]}) := \max_{\indexProsumer \in [\setProsumer]} \; \frac{
    \objectiveProsumer{\indexProsumer}(\decisionProsumer{\indexProsumer})}{\valueProsumer{\indexProsumer}},\label{obj:minimax2}
\end{align}
\end{subequations}
\new{which we refer to as the \textit{Scaled Minimax approach}.
Note that in both the scaled minimax approach ($\objectiveMinimaxSavings{\setProsumer}$) and the proportional approach ($\objectiveProportional{\setProsumer}$), there are multiple solutions with different aggregated costs.}
\new{We assume here we have defined a way to select a solution among them. }

\subsection{Acceptability constraints}
\label{subsec:acceptibility_deterministic}

\change{Having delineated several methodologies for equitable cost distribution,}
we must convince \change{agents} to be part of the aggregation. 
We consider that \change{agents are individually rational, that is a contract cannot be deemed acceptable if at least one agent would} be better off independently \new{\ie $ \valueProsumer{\indexProsumer} \leq \objectiveProsumer{\indexProsumer}(\decisionProsumer{\indexProsumer})$, where $\objectiveProsumer{\indexProsumer}(\decisionProsumer{\indexProsumer})$ is the cost of $\indexProsumer$ in the aggregation.} 
We can go one step further and require that, to find the contract acceptable, they benefit from it\new{, \ie $\valueProsumer{\indexProsumer} > \objectiveProsumer{\indexProsumer}(\decisionProsumer{\indexProsumer})$}. 
We thus define the acceptability set $\acceptableSet{\coefficientAcceptability}{\indexProsumer}$ appearing in~\eqref{eq:abstract_acceptability_constraint} as follows:
\begin{equation}
    \acceptableSet{\coefficientAcceptability}{\indexProsumer} := \ba{\; \utilityProsumer{\indexProsumer} \, | \;  \utilityProsumer{\indexProsumer} \leq \coefficientAcceptability \; \valueProsumer{\indexProsumer} \;} ,
    \label{eq:acceptability}
\end{equation}
where $\coefficientAcceptability  \in (0,1]$ is given.
\new{Then, we say a solution is $\coefficientAcceptability-$acceptable if it is contained in $\acceptableSet{\coefficientAcceptability}{\indexProsumer}$.}
%
%
Acceptability sets are independent from one \change{agent} to another. 
\change{
We then define \textit{global acceptability} as the cartesian product of all acceptability sets 
$\acceptableSet{\coefficientAcceptability}{ }:=\acceptableSet{\coefficientAcceptability}{\indexProsumer_1}\times \dots \times \acceptableSet{\coefficientAcceptability}{\indexProsumer_{\setProsumer}}$.
}

\new{
Enforcing acceptability constraints to $\modelAggregator$ restricts the feasible solution set, potentially leading to higher aggregated cost. 
We define the \textit{price of acceptability} as 
\begin{equation}
    PoA:=v^\star_{\acceptableSet{\coefficientAcceptability}{ }} - v_\emptyset^\star,
\end{equation} where $v^\star_{\acceptableSet{\coefficientAcceptability}{ }}$ is the optimal value of $\modelAggregator$ with acceptability constraints $\acceptableSet{\coefficientAcceptability}{ }$, and $v^\star$ is the optimal value of $\modelAggregator$ without them. 
}

\vspace{2mm}
\begin{remark}
    In the scaled minimax model with \operatorName $\objectiveMinimaxSavings{\setProsumer}$, the optimal solution is $1-$acceptable. 
    Indeed, if the agents don't take advantage of the aggregation, then $\objectiveProsumer{\indexProsumer}(\decisionProsumer{\indexProsumer})=\valueProsumer{\indexProsumer}$ and we get a \new{feasible solution with respect to $\acceptableSet{\coefficientAcceptability}{ }$} of optimal value $1$. 
    Further, if we consider the problem: 
\begin{subequations}
\begin{align}
    \Min_{\alpha} \qquad & \; \coefficientAcceptability \\
    \text{s.t.} \qquad & \decisionProsumer{\indexProsumer} \in \constraintProsumer{\indexProsumer} & \forall \indexProsumer \in [\setProsumer] \\
    &  \change{h(\decisionProsumer{1}, \dots, \decisionProsumer{\setProsumer})} \in \constraintMarket  \\
    & \objectiveProsumer{\indexProsumer}(\decisionProsumer{\indexProsumer}) \in \acceptableSet{\coefficientAcceptability}{\indexProsumer} & \forall \indexProsumer \in [\setProsumer],
\end{align}
\label{eq:problem_alpha}
\end{subequations}
it is equivalent to problem $\modelAggregator$ with \operatorName $\objectiveMinimaxSavings{\setProsumer}$ \new{with no acceptability constraints:
\begin{subequations}
\begin{align}
    \Min_{x} \qquad & \; \Max_{\indexProsumer \in [\setProsumer]} \qquad \frac{\objectiveProsumer{\indexProsumer}(\decisionProsumer{\indexProsumer})}{\valueProsumer{\indexProsumer}} \label{eq:proportional_maximum}\\
    \text{s.t.} \qquad & \decisionProsumer{\indexProsumer} \in \constraintProsumer{\indexProsumer} & \forall \indexProsumer \in [\setProsumer] \\
    &  h(\decisionProsumer{1}, \dots, \decisionProsumer{\setProsumer}) \in \constraintMarket.
\end{align}
\end{subequations}
Then, when we linearize the maximum in \eqref{eq:proportional_maximum}, we fall back into Problem~\eqref{eq:problem_alpha}, by definition of $\acceptableSet{\coefficientAcceptability}{\indexProsumer}$.
}
\label{remark:acceptable_minimax}
\end{remark}

\vspace{2mm}
\begin{remark}
    Note that our problem with the proportional operator $\objectiveProportional{\setProsumer}$ necessarily yields a solution $(1-\epsilon)-$acceptable, with $\epsilon>0$. 
    Indeed, if for agent $\indexProsumer$, $\objectiveProsumer{\indexProsumer}(\decisionProsumer{\indexProsumer})\geq \valueProsumer{\indexProsumer}$, then $\log(\valueProsumer{\indexProsumer}-\objectiveProsumer{\indexProsumer}(\decisionProsumer{\indexProsumer}))$ is undefined. 
\label{remark:acceptable_proportional}
\end{remark}

\vspace{2mm}
For simplicity, in the rest of the paper, we assume $\coefficientAcceptability=1$.
We later discuss how to extend the acceptability constraint to a dynamic (see \cref{subsec:dynamic_acceptability}) and stochastic framework (see \cref{subsec:stochastic_acceptability}).
Finally, 
combining different objective functions with acceptability constraints, we observe on a small illustration their impact on the solution in the following section.

\section{Application to consumer aggregation on the day-ahead and balancing market}
\label{subsec:deterministic_illustration}
\new{In this section, we adapt and illustrate the framework presented in \cref{sec:problem} to the problem of prosumers aggregation on electricity markets.
More specifically, the prosumers have access to:
the \textit{day-ahead market}, where every day at 2 pm, prices and electrical energies are set for all across Europe for the twenty-four hours of the next day;
and the \textit{balancing market} on which prosumers must buy or sell electricity at real-time prices to ensure power system balance.
A minimum trade of 0.1 MWh of energy is required to participate in the day-ahead market. 
}

\begin{table}[ht]
\centering
\begin{subtable}{.4\linewidth}
    \renewcommand{\arraystretch}{1.5}
    \centering
    \begin{tabular}{|c||*{5}{c|}}\hline
        $\indexStage$ & \textbf{1} & \textbf{2} & \textbf{3} & \textbf{4} & \textbf{5}\\ 
        \hline\hline
        $\pricesDAExample{\indexStage}$  & 2 & 16 & 1 & 10 & 1  \\\hline
        $\pricesBExample{\indexStage}$ & 6 & 25 & 5 & 15 & 5  \\\hline
        $\minDA{\indexStage}$ & 11 & 11 & 11 & 11 & 11 \\\hline
    \end{tabular}

\end{subtable}%
\hspace{0.3cm}
\begin{subtable}{.4\linewidth}
    \centering
    \def\arraystretch{1.5}

    \begin{tabular}{|c||*{4}{c|}}\hline
        
         &  $A_1$ & $A_2$ & $A_3$ & $A_4$ \\\hline\hline
        $\quantityMinEnergyExample{\indexProsumer}$ & 0 & 5 & 0 & 2 \\\hline 
        $\quantityMaxEnergyExample{\indexProsumer}$ & 5 & 5 & 4 & 3 \\\hline
        $\quantityTotalEnergyExample{\indexProsumer}$ & 10 & 25 & 8 & 15 \\ 
        \hline
    \end{tabular}
\end{subtable}
\caption{Parameters values}
\label{tab:parameters}
\end{table}

\change{We propose a toy model to illustrate} the implications of each model proposed in \cref{sec:problem}.
Therefore, we consider a problem with four consumers ($\setProsumer=4$) on five stages ($\nbStages=5$).
At each stage $\indexStage$, we must decide how much energy $\quantityEnergyDAExample{\indexStage}{\indexProsumer}$ (\textit{resp.} $\quantityEnergyBExample{\indexStage}{\indexProsumer}$) to purchase from the day-ahead (\textit{resp.} balancing) market for consumer $\indexProsumer$
\new{. Thus in $\modelProsumer{\indexProsumer}$, we have $\decisionProsumer{\indexProsumer}:=(\quantityEnergyDAExample{\indexStage}{\indexProsumer},\quantityEnergyBExample{\indexStage}{\indexProsumer})$}.
Each consumer has bounds $\nc{\quantityMinEnergyExample{\indexProsumer} ; \quantityMaxEnergyExample{\indexProsumer}}$ on its electricity consumption, and a total consumption $\quantityTotalEnergyExample{\indexProsumer}$ to meet at the end of the horizon, \new{amounting to feasible set $\constraintProsumer{\indexProsumer}$}.
\new{Note that the upper bounds on electricity consumption simplify physical constraints that would ensure a finite volume of traded electricity.}
We introduce binary variables $\buyDA{\indexStage}$, representing the decision to buy a day-ahead, to model the minimum volume requirement for the day-ahead market, \new{which composes the external constraints $\cM$}.
The objective for consumer $\indexProsumer$ is to minimize its \change{electricity} costs:
\begin{subequations}
\begin{align}
    \objectiveProsumer{\indexProsumer}(\decisionProsumer{\indexProsumer}) = \sum_{\indexStage=1}^\nbStages \nc{\, \pricesDAExample{\indexStage} \quantityEnergyDAExample{\indexStage}{\indexProsumer} + \pricesBExample{\indexStage} \quantityEnergyBExample{\indexStage}{\indexProsumer} \,} ,
\end{align}
where $\pricesDAExample{\indexStage}$ (\textit{resp.} $\pricesBExample{\indexStage}$) is the price of \change{electricity} at $\indexStage$ on the day-ahead (\textit{resp.} balancing) market.
We obtain the simple prosumer $\modelProsumer{\indexProsumer}$ and the aggregated model $\modelAggregator$: 

\begin{alignat}{5}
    \hspace{-2.cm} \modelProsumer{\indexProsumer} \; \Min_{\decisionProsumer{\indexProsumer}} \quad & \objectiveProsumer{\indexProsumer}(\decisionProsumer{\indexProsumer}) & 
    \hspace{-0.4cm} \modelAggregator \; \Min_{\decisionProsumer{ }} \; & \operatorAggregation{\setProsumer}\Bp{(\objectiveProsumer{\indexProsumer}(\decisionProsumer{\indexProsumer}))_{\indexProsumer \in [\setProsumer]}} \\
    \text{s.t.} \quad &  \quantityMinEnergyExample{\indexProsumer} \leq \quantityEnergyDAExample{\indexStage}{\indexProsumer} + \quantityEnergyBExample{\indexStage}{\indexProsumer} \leq  \quantityMaxEnergyExample{\indexProsumer} \quad \forall \indexStage \qquad \quad & &  
    \text{s.t.} \quad 
    \quantityMinEnergyExample{\indexProsumer} \leq \quantityEnergyDAExample{\indexStage}{\indexProsumer} + \quantityEnergyBExample{\indexStage}{\indexProsumer} \leq  \quantityMaxEnergyExample{\indexProsumer} \quad \forall \indexStage \; \forall \indexProsumer \\
    & \sum_{\indexStage=1}^\nbStages (\quantityEnergyDAExample{\indexStage}{\indexProsumer} + \quantityEnergyBExample{\indexStage}{\indexProsumer}) \geq \quantityTotalEnergyExample{\indexProsumer}  & & \qquad \sum_{\indexStage=1}^\nbStages (\quantityEnergyDAExample{\indexStage}{\indexProsumer} + \quantityEnergyBExample{\indexStage}{\indexProsumer}) \geq \quantityTotalEnergyExample{\indexProsumer} \quad \forall \indexProsumer  \\
    & \minDA{\indexStage} \buyDA{\indexStage} \leq\quantityEnergyDAExample{\indexStage}{\indexProsumer}  \leq M \, \buyDA{\indexStage} \quad \forall \indexStage
    & & \quad \minDA{\indexStage} \buyDA{\indexStage} \leq \sum_{\indexProsumer \in [\setProsumer]}  \quantityEnergyDAExample{\indexStage}{\indexProsumer}  \leq M \, \buyDA{\indexStage} \; \forall \indexStage \label{eq:minimum_requirement} \\
    & \buyDA{\indexStage} \in \na{0,1} \quad \forall \indexStage,
    & & \qquad \buyDA{\indexStage} \in \na{0,1} \quad \forall \indexStage, 
    \label{eq:illustration_deterministic}
\end{alignat}
\end{subequations}

where $\operatorAggregation{ }$ is the chosen \operatorName for the aggregation. 
We solve this small problem with the utilitarian operator $\objectiveUtilitarian{\setProsumer}$, with the scaled minimax operator $\objectiveMinimaxSavings{\setProsumer}$ and with the proportional operator $\objectiveProportional{\setProsumer}$.
For all \operatorName, we solve the problem with and without acceptability constraints, with $\coefficientAcceptability=1$.
\new{In the proportional and scaled minimax approaches, we do not optimize the aggregated costs.
As the optimal solution is not necessarily unique, there can be different optimal solutions with different aggregated costs.
Among those, we choose one with minimum aggregated costs.}

\new{We refer to the model with agent operator $f \in \na{U, SM, P}$ (for Utilitarian, Scaled Minimax and Proportional) and acceptability constraints set $a \in \na{\emptyset, \coefficientAcceptability}$ (for no acceptability constraints, or acceptability constraints given by $\cA_\coefficientAcceptability$) as $\modeName{f}{a}$, and $\cA=\emptyset$ corresponds to a model without acceptability constraints.}
\new{Finally, we compute here Shapley's values (see \cref{ann:shapley} for more details), which are commonly recognized as a fair solution, to compare them to the solutions we obtain with our models. }

We show on a small artificial illustration how all these models can lead to different solutions. 
\new{
Each model can be evaluated through two metrics: 
first, the efficiency of the model \ie the overall costs of the aggregation;
second, the fairness of the model \ie how distributed are the costs over prosumers. 
}
For the prosumers' parameters and market prices we use the data on \cref{tab:parameters}.
\change{We observe the allocation of costs over consumers on \cref{fig:results_deterministic_obj}, the resulting percentage of savings made by each consumer on \cref{tab:results_deterministic}, and the detail of day-ahead and balancing purchases on \cref{tab:details_purchases}.}

First, it's worth noting that none of the consumers can individually access the day-ahead market \new{as for any prosumer $\quantityMaxEnergyExample{\indexProsumer} \leq \minDA{\indexStage}$ and thus constraint~\eqref{eq:minimum_requirement} excludes any purchase on the day-ahead market.}
In the utilitarian model $\modeName{U}{\emptyset}$,the primary focus is to minimize aggregated costs,
making it optimal to always consistently access the day-ahead market as a group. 
To achieve this, consumer $\consumer{1}$ redistributes its energy load across $4$ time steps, incurring a higher individual cost ($64\%$ higher) than when acting independently. 
By adding acceptability constraints to the model ($\modeName{U}{\coefficientAcceptability}$), \new{the model loses in efficiency but now satisfies individual rationality:
$PoA=3\%$ of $\modeName{U}{\emptyset}$'s costs. 
}
We observe that the aggregated costs of consumers \change{slightly increases}, 
but now the charge of energy needed to access the day-ahead market is shared between $\consumer{1}$ and $\consumer{3}$,
although $\consumer{3}$ does not gain anything from the the aggregation ($0\%$ of savings). 

\begin{figure}[h]
\label{table:details_purchases}
    \centering
    \scalebox{1.2}{\pgfplotsset{compat=1.9}

\begin{tikzpicture}

\pgfplotsset{
        show sum on top/.style={
            /pgfplots/scatter/@post marker code/.append code={%
                \node[
                    at={(normalized axis cs:%
                            \pgfkeysvalueof{/data point/x},%
                            \pgfkeysvalueof{/data point/y})%
                    },
                    anchor=south,
                ]
                {\pgfmathprintnumber{\pgfkeysvalueof{/data point/y}}};
            },
            font=\footnotesize
        },
    }
    
 \begin{axis}[
     ybar stacked,
     bar width=0.5cm,
     ylabel={Total Costs (\$)},
     xtick={1, 4, 6, 9, 11, 14, 16},
     xticklabels={{\scriptsize $\modelProsumer{\indexProsumer}$, 
     \scriptsize $\modeName{U}{\emptyset}$, 
     \scriptsize $\modeName{U}{\coefficientAcceptability}$, 
     \scriptsize $\modeName{SM}{\emptyset}$, 
     \scriptsize $\modeName{SM}{\coefficientAcceptability}$, 
     \scriptsize $\modeName{P}{\emptyset}$,
     \scriptsize $\modeName{P}{\coefficientAcceptability}$}},
    legend style={at={(0.5,0.95)},
        anchor=north,
        legend columns=4,
        draw=none
     },
 ]

 \addplot +[draw=blue!70,fill=blue!40, nodes near coords style={blue, font=\footnotesize, anchor=center, yshift=0.1cm}, nodes near coords, nodes near coords align={center}, point meta=explicit symbolic] coordinates {(1,50.0) (4,82.0) (6,26.0) (9,20.1) (11,20.1) (14,13.0) (16,13.0) };
        
\addplot +[pattern = crosshatch, pattern color =red!50,draw=red!60, nodes near coords align={center}, nodes near coords style={red, font=\footnotesize, anchor=center, yshift=0.1cm}, nodes near coords,point meta=explicit symbolic] coordinates {(1,280.0) (4,150.0) (6,175.0) (9,195.0) (11,195.0) (14,220.0) (16,220.0)};

\addplot +[fill=brown!20,draw=brown!70, nodes near coords align={center}, nodes near coords style={brown, font=\footnotesize, anchor=center, yshift=0.1cm}, nodes near coords,point meta=explicit symbolic] coordinates {(1,40.0) (4,11.0) (6,40.0) (9,27.9) (11,27.9) (14,8.0) (16,8.0)};

\addplot +[pattern = north west lines, pattern color=orange!70,draw=orange!70, nodes near coords align={center}, nodes near coords style={orange, font=\footnotesize, anchor=center, yshift=0.1cm},nodes near coords,point meta=explicit symbolic, show sum on top] coordinates {(1,168.0) (4,90.0) (6,105.0) (9,117.0) (11,117.0) (14,132.0) (16,132.0)};

\legend{$\consumer{1}$, $\consumer{2}$, $\consumer{3}$, $\consumer{4}$}
\end{axis}

\end{tikzpicture}}
    \caption{We observe the result of the static Problem~\eqref{eq:illustration_deterministic} with parameters given in \cref{tab:parameters}. 
    The \change{bars} correspond to the outcome of different models, the number above being the total cost.
    The first bar is the non-aggregated model: we solve each $\modelProsumer{\indexProsumer}$ independently. 
    Then, there are three groups of two \change{bars}, each group corresponding to a choice of \operatorName ($\objectiveUtilitarian{\setProsumer},\objectiveMinimaxSavings,\objectiveProportional{\setProsumer}$).
    Then, \change{for each objective function, we present the results of the model}, first without and then with, acceptability constraints $\acceptableSet{\coefficientAcceptability}{ }$,  with $\coefficientAcceptability=1$.  
    Each \change{bar} is decomposed in $4$ blocks corresponding to the cost incurred by each consumer $i$.
    \new{At the top of each bar, we can read the sum of aggregated costs in the corresponding model.}
    }
    \label{fig:results_deterministic_obj}
\end{figure}

\begin{table}[!ht]
\caption{Percentage of savings $\frac{\valueProsumer{\indexProsumer} - \objectiveProsumer{\indexProsumer}(\decisionProsumer{\indexProsumer})}{\valueProsumer{\indexProsumer}}$ made by $A_\indexProsumer$ in \new{the model $\modeName{f}{a}$ depending on agent operator $f \in \na{U,SM,P}$ and acceptability set $a \in \na{\emptyset, \coefficientAcceptability}$} and $PoA$ (in percentage) of the corresponding model.}
\label{tab:results_deterministic}
\centering
\renewcommand{\arraystretch}{1.2}
\setlength{\tabcolsep}{4pt}
\begin{tabular}{c|ccccc|ccccc|ccccc}
\toprule
 & \multicolumn{5}{c}{\footnotesize Utilitarian \new{$\objectiveUtilitarian{\setProsumer}$}} & \multicolumn{5}{c}{\footnotesize Minimax \new{$\objectiveMinimaxSavings{\setProsumer}$}}  & \multicolumn{5}{c}{\footnotesize Proportional \new{$\objectiveProportional{\setProsumer}$}}  \\
\cmidrule(lr){2-6} \cmidrule(lr){7-11} \cmidrule(lr){12-16}
 & A1 & A2 & A3 & A4 & $PoA$ & A1 & A2 & A3 & A4 & $PoA$ & A1 & A2 & A3 & A4 & $PoA$ \\
\midrule
\new{$\emptyset$} & -64 & 46 & 72 & 46 & 0 & 60 & 30 & 30 & 30 & 0 & 74 & 21 & 80 & 21 & 0  \\
\new{$\acceptableSet{\coefficientAcceptability}{ }$}  & 48 & 37 & 0 & 37 & 4 & 60 & 30 & 30 & 30 & 0 & 74 & 21 & 80 & 21 & 0 \\ \hline \hline 
Shapley & 114 & 20 & 111 & 28 & 0 & & & & & & & & & & \\
\bottomrule
\end{tabular}
\end{table}

\new{
\begin{table}[ht]
\renewcommand{\arraystretch}{1.1}
\setlength{\tabcolsep}{4pt}

\begin{subtable}{.45\linewidth}
\centering
\begin{tabular}{l|cc|cc|cc|cc}
        \hline
         & \multicolumn{2}{c}{A1} & \multicolumn{2}{c}{A2} & \multicolumn{2}{c}{A3} & \multicolumn{2}{c}{A4} \\ \hline
        t & DA & B & DA & B &  DA & B & DA & B \\
1 & \textcolor{blue}{\textit{\textbf{  0 }}} &  \; 0 \; &  \; 5 \; &  \; 0 \; &  \; 3 \; &  \; 0 \; &  \; 3 \; &  \; 0 \; \\
2 &  \; 3 \; &  \; 0 \; &  \; 5 \; &  \; 0 \; & \textcolor{blue}{\textit{\textbf{  0 }}} &  \; 0 \; &  \; 3 \; &  \; 0 \; \\
3 &  \; 2 \; &  \; 0 \; &  \; 5 \; &  \; 0 \; & \; 1 \; &  \; 0 \; &  \; 3 \; &  \; 0 \; \\
4 &  \; 3 \; &  \; 0 \; &  \; 5 \; &  \; 0 \; & \textcolor{blue}{\textit{\textbf{  0 }}} &  \; 0 \; &  \; 3 \; &  \; 0 \; \\
5 &  \; 2 \; &  \; 0 \; &  \; 5 \; &  \; 0 \; & \; 4 \; &  \; 0 \; &  \; 3 \; &  \; 0 \; \\
\hline
\end{tabular}
\caption{$\modeName{U}{\emptyset}$ \textcolor{white}{$\modeName{U}{\coefficientAcceptability}$}}
\end{subtable}%
\begin{subtable}{.45\linewidth}
\centering
\begin{tabular}{||cc|cc|cc|cc}
    \hline
    
    \multicolumn{2}{c}{A1} & \multicolumn{2}{c}{A2} & \multicolumn{2}{c}{A3} & \multicolumn{2}{c}{A4} \\  \hline
    DA & B & DA & B &  DA & B &  DA & B \\
\textcolor{blue}{\textit{\textbf{  0 }}} &  \; 0 \; &  \; 5 \; &  \; 0 \; &  \; 3 \; &  \; 0 \; &  \; 3 \; &  \; 0 \; \\
1.01 &  \; 0 \; &  \; 5 \; &  \; 0 \; & 1.93 &  \; 0 \; &  \; 3 \; &  \; 0 \; \\
3.93 &  \; 0 \; &  \; 5 \; &  \; 0 \; & \textcolor{blue}{\textit{\textbf{  0 }}} &  \; 0 \; &  \; 3 \; &  \; 0 \; \\
\textcolor{blue}{\textit{\textbf{  0 }}} &  \; 0 \; &  \; 0 \; & \textcolor{red}{\textit{\textbf{  5 }}} & \textcolor{blue}{\textit{\textbf{  0 }}} &  \; 0 \; &  \; 0 \; & \textcolor{red}{\textit{\textbf{  3 }}} \\
 \; 5 \; &  \; 0 \; &  \; 5 \; &  \; 0 \; &  \; 3 \; &  \; 0 \; &  \; 3 \; &  \; 0 \; \\
\hline \end{tabular}
\caption{$\modeName{U}{\coefficientAcceptability}$}
\end{subtable}

\begin{subtable}{.45\linewidth}
\centering
    \begin{tabular}{l|cc|cc|cc|cc}
    \hline  & \multicolumn{2}{c}{A1} & \multicolumn{2}{c}{A2} & \multicolumn{2}{c}{A3} & \multicolumn{2}{c}{A4} \\  \hline
    t & DA & B & DA & B &  DA & B &  DA & B \\
1 & \textcolor{blue}{\textit{\textbf{  0 }}} &  \; 0 \; &  \; 5 \; &  \; 0 \; &  \; 3 \; &  \; 0 \; &  \; 3 \; &  \; 0 \; \\
2 & \textcolor{blue}{\textit{\textbf{  0 }}} &  \; 0 \; &  \; 0 \; & \textcolor{red}{\textit{\textbf{ 5 }}} & \textcolor{blue}{\textit{\textbf{  0 }}} &  \; 0 \; &  \; 0 \; & \textcolor{red}{\textit{\textbf{ 3 }}} \\
3 & 4.44 &  \; 0 \; &  \; 5 \; &  \; 0 \; & 1.56 &  \; 0 \; &  \; 3 \; &  \; 0 \; \\
4 & 1.13 &  \; 0 \; &  \; 5 \; &  \; 0 \; & 1.87 &  \; 0 \; &  \; 3 \; &  \; 0 \; \\
5 & 4.44 &  \; 0 \; &  \; 5 \; &  \; 0 \; & 1.56 &  \; 0 \; &  \; 3 \; &  \; 0 \; \\
\hline \end{tabular}
\caption{$\modeName{SM}{\emptyset}$ \textcolor{white}{$\modeName{U}{\coefficientAcceptability}$}}
\end{subtable}%
\begin{subtable}{.45\linewidth}
        \centering
        \begin{tabular}{||cc|cc|cc|cc}
        \hline  \multicolumn{2}{c}{A1} & \multicolumn{2}{c}{A2} & \multicolumn{2}{c}{A3} & \multicolumn{2}{c}{A4} \\  \hline
        DA & B & DA & B &  DA & B &  DA & B \\
\textcolor{blue}{\textit{\textbf{  0 }}} &  \; 0 \; &  \; 5 \; &  \; 0 \; &  \; 3 \; &  \; 0 \; &  \; 3 \; &  \; 0 \; \\
\textcolor{blue}{\textit{\textbf{  0 }}} &  \; 0 \; &  \; 0 \; & \textcolor{red}{\textit{\textbf{  5 }}} & \textcolor{blue}{\textit{\textbf{  0 }}} &  \; 0 \; &  \; 0 \; & \textcolor{red}{\textit{\textbf{  3 }}} \\
4.44 &  \; 0 \; &  \; 5 \; &  \; 0 \; & 1.56 &  \; 0 \; &  \; 3 \; &  \; 0 \; \\
1.13 &  \; 0 \; &  \; 5 \; &  \; 0 \; & 1.87 &  \; 0 \; &  \; 3 \; &  \; 0 \; \\
4.44 &  \; 0 \; &  \; 5 \; &  \; 0 \; & 1.56 &  \; 0 \; &  \; 3 \; &  \; 0 \; \\
\hline \end{tabular}
\caption{$\modeName{SM}{\coefficientAcceptability}$}
\end{subtable}

\begin{subtable}{.45\linewidth}
    \centering
    \begin{tabular}{l|cc|cc|cc|cc}
    \hline  & \multicolumn{2}{c}{A1} & \multicolumn{2}{c}{A2} & \multicolumn{2}{c}{A3} & \multicolumn{2}{c}{A4} \\ \hline
    t & DA & B & DA & B &  DA & B &  DA & B \\
1 &  \; 3 \; &  \; 0 \; &  \; 5 \; &  \; 0 \; & \textcolor{blue}{\textit{\textbf{  0 }}} &  \; 0 \; &  \; 3 \; &  \; 0 \; \\
2 & \textcolor{blue}{\textit{\textbf{  0 }}} &  \; 0 \; &  \; 0 \; & \textcolor{red}{\textit{\textbf{  5 }}} & \textcolor{blue}{\textit{\textbf{  0 }}} &  \; 0 \; &  \; 0 \; & \textcolor{red}{\textit{\textbf{  3 }}} \\
3 & 3.5 \; &  \; 0 \; &  \; 5 \; &  \; 0 \; & \; 4 \; &  \; 0 \; &  \; 3 \; &  \; 0 \; \\
4 & \textcolor{blue}{\textit{\textbf{  0 }}} &  \; 0 \; &  \; 0 \; & \textcolor{red}{\textit{\textbf{  5 }}} & \textcolor{blue}{\textit{\textbf{  0 }}} &  \; 0 \; &  \; 0 \; & \textcolor{red}{\textit{\textbf{  3 }}} \\
5 & 3.5 \; &  \; 0 \; &  \; 5 \; &  \; 0 \; & \; 4 \; &  \; 0 \; &  \; 3 \; &  \; 0 \; \\
\hline \end{tabular}
\caption{$\modeName{P}{\emptyset}$ \textcolor{white}{$\modeName{U}{\coefficientAcceptability}$}}
\end{subtable}%
\begin{subtable}{.45\linewidth}
    \centering
    \begin{tabular}{||cc|cc|cc|cc}
    \hline  \multicolumn{2}{c}{A1} & \multicolumn{2}{c}{A2} & \multicolumn{2}{c}{A3} & \multicolumn{2}{c}{A4} \\ \hline
    DA & B & DA & B &  DA & B &  DA & B \\
 \; 3 \; &  \; 0 \; &  \; 5 \; &  \; 0 \; & \textcolor{blue}{\textit{\textbf{  0 }}} &  \; 0 \; &  \; 3 \; &  \; 0 \; \\
\textcolor{blue}{\textit{\textbf{  0 }}} &  \; 0 \; &  \; 0 \; & \textcolor{red}{\textit{\textbf{  5 }}} & \textcolor{blue}{\textit{\textbf{  0 }}} &  \; 0 \; &  \; 0 \; & \textcolor{red}{\textit{\textbf{  3 }}} \\
3.5 \; &  \; 0 \; &  \; 5 \; &  \; 0 \; & \; 4 \; &  \; 0 \; &  \; 3 \; &  \; 0 \; \\
\textcolor{blue}{\textit{\textbf{  0 }}} &  \; 0 \; &  \; 0 \; & \textcolor{red}{\textit{\textbf{  5 }}} & \textcolor{blue}{\textit{\textbf{  0 }}} &  \; 0 \; &  \; 0 \; & \textcolor{red}{\textit{\textbf{  3 }}} \\
3.5 \; &  \; 0 \; &  \; 5 \; &  \; 0 \; & \; 4 \; &  \; 0 \; &  \; 3 \; &  \; 0 \; \\
\hline \end{tabular}
\caption{$\modeName{P}{\coefficientAcceptability}$}
\end{subtable}
\caption{Quantity of energy purchased on the day-ahead and balancing markets per stage for all consumers depending on different models: in bold italic, we highlight purchases on the (more expensive) balancing market and stages where no purchases are made on the day-ahead market.}
\label{tab:details_purchases}
\end{table}

}

\change{Conversely, the proportional solution (from model $\modeName{P}{\emptyset}$) adopts a more bargaining-oriented approach, resulting in collaboration only during time slots ($t\in \na{1, 3, 5}$) with lower prices. 
Indeed, as $\consumer{1}$ and $\consumer{3}$ are not forced to consume energy at all times ($\quantityMinEnergyExample{1}=\quantityMinEnergyExample{3}=0$), they can shift their consumption to time slots with lower prices. 
On the contrary, $\consumer{2}$ and $\consumer{4}$ must always consume energy ($\quantityMinEnergyExample{2}=5, \quantityMinEnergyExample{4}=2$), and the two of them together cannot access the day-ahead market either.
Thus, in $\modeName{P}{\emptyset}$, the solution is for $\consumer{1}$ and $\consumer{3}$ to consume only in time steps $\na{1,3,5}$, 
which leaves $\consumer{2}$ and $\consumer{4}$ to operate independently at $t=2, t=4$, resulting in limited savings ($21\%$) compared to the scaled minimax approach ($\modeName{SM}{\emptyset}$).}
As noticed in  \cref{remark:acceptable_proportional}, the solution is necessarily 1-acceptable.
Therefore, the solution is the same in $\modeName{P}{\emptyset}$ and $\modeName{P}{\coefficientAcceptability}$.
Moreover, the proportional solution yields the worst aggregated costs \ie the less efficient solution. 

\change{With the scaled minimax approach, the model $\modeName{SM}{\emptyset}$ yields a trade-off between efficiency and fairness compared to $\modeName{U}{\coefficientAcceptability}$:
we observe that $\consumer{1}$ and $\consumer{3}$ decide to stop consuming at expensive time steps, thus achieving greater savings.
The model also encourages more cooperation than the proportional model $\modeName{P}{\emptyset}$, as we can observe on \cref{tab:details_purchases}.
As a result, in this model, all consumers achieve similar proportional savings, amounting to approximately 30\% compared to operating independently, 
at the exception of $\consumer{1}$ that can save up to $60\%$. 
This means that any solution where $\consumer{1}$ shifts its consumption to other time slots to help others access the day-ahead market,
 would increase its costs too much, and $\consumer{3}$ would save less than $30\%$: this cannot be an optimal solution of $\modeName{SM}{\emptyset}$. 
However, the aggregated cost of the aggregation is higher than with $\modeName{U}{\emptyset}$ and $\modeName{U}{\coefficientAcceptability}$.
}
Again, adding acceptability constraints does not change the solution, as the scaled minimax problem is innately 1-acceptable (see \Cref{remark:acceptable_minimax}).

\new{
Lastly, we observe on \cref{tab:results_deterministic} the allocation of savings with a post-allocation rule based on Shapley's values. 
This approach leverages the efficiency of $\modeName{U}{\emptyset}$ but then re-allocates costs to obtain a fair and acceptable solution.
In this application, $\consumer{1}$ and $\consumer{3}$ save respectively $114\%$ and $111\%$ of their costs compared to operating independently, 
which amounts to them being paid by the aggregation to participate.
Even though this allows $\consumer{2}$ and $\consumer{4}$ to gain from the aggregation, this questions the acceptability of this solution as some would earn money while others have residual costs. 
Furthermore, this method becomes impractical when dealing with large problems involving many prosumers and time steps due to the extensive computational requirements.
Additionally, adapting this method to dynamic and stochastic contexts is unclear, so we do not consider it further.

}

\section{Fairness \change{accross time}}
\label{sec:dynamic}

In most use cases, we can assume that the aggregation of \change{agents} is thought to stay in place over long periods. 
One of the challenges of this long-term setting is incentivizing \change{agents} not to leave the aggregation, which requires adjusting the acceptability constraints of the static case.

\subsection{Problem formulation}
\label{subsec:dynamic_problem}

We consider a problem with $\horizon$ stages corresponding to consecutive times where decisions are made. 
At each stage $\indexTime \in [\horizon]$, \change{agent} $\indexProsumer$ makes a decision $\decisionProsumerDynamic{\indexTime}{\indexProsumer} \in \constraintProsumerDynamic{\indexTime}{\indexProsumer}$, incurring a cost 
$\instantaneousCostProsumer{\indexTime}{\indexProsumer}(\decisionProsumerDynamic{\indexTime}{\indexProsumer})$.
Those stage costs are aggregated through a time operator $\operatorTime{\horizon}{\indexProsumer}:\RR^T \to \RR$.
Thus, the \change{agent} $\indexProsumer$'s problem reads:
\begin{subequations}
    \begin{align}
        \modelProsumerDynamic{\horizon}{\indexProsumer} \; := \quad \Min_{\decisionProsumerDynamic{\indexTime}{\indexProsumer}} \quad & 
        \operatorTime{\horizon}{\indexProsumer} \Bp{ (\instantaneousCostProsumer{\indexTime}{\indexProsumer}(\decisionProsumerDynamic{\indexTime}{\indexProsumer}))_{\indexTime \in [\horizon]}
        }\label{dynamic-obj} \\
        \text{s.t} \quad & \decisionProsumerDynamic{\indexTime}{\indexProsumer} \in \constraintProsumerDynamic{\indexTime}{\indexProsumer} & \forall \indexTime \label{feasible-space2} \\
        &\change{\decisionProsumerDynamic{\indexTime}{\indexProsumer}} \in \constraintMarketDynamic{\indexTime} & \forall \indexTime .
    \end{align}
    \label{pb:prosumer_dynamic}
\end{subequations}

A typical example of time-aggregator $\operatorTime{\horizon}{\indexProsumer}$ is the \change{discounted} sum of stage costs \ie  dropping the dependence in $x_\indexProsumer$ for clarity's sake:
 $$ \operatorTime{\horizon}{\indexProsumer}((L_t)_{t \in [\horizon]})=\sum_{\indexTime \in [\horizon]} r^\indexTime L_\indexTime,$$
 for $r\in(0,1]$.
Alternatively, $\operatorTime{\horizon}{\indexProsumer}$ can be defined as the maximum of stage costs. 
This might happen for electricity markets where a prosumer aims at \emph{peak shaving} \ie minimizing peak electricity demand.
Further, time-aggregation operators may vary among agents, who may express different sensitivity to time \new{ \ie the discounted rate $r$ varies among prosumers.}
\delete{In the remainder of the paper, for simplicity, we always consider the time-aggregator as the sum of stage-costs, where discounted rates are directly included in the definition of the stage-cost $\instantaneousCostProsumer{\indexTime}{\indexProsumer}$.}


We now write the aggregation problem within this framework.
Note that we can cast the current multistage setting into the setting of \cref{sec:problem},
by \change{decomposing each agent into $\horizon$ independent stage-wise sub-agents: 
then we have $\setProsumer \times \horizon$ and we can use the methodology of \cref{sec:problem}.}
Thus we need to define an operator  $\operatorAggregation{\setProsumer \times \horizon}$ that takes $\na{\instantaneousCostProsumer{\indexTime}{\indexProsumer}}_{\indexTime \in [\horizon],\indexProsumer \in [\setProsumer]}$ as input. 

However, in most settings, it is reasonable to assume that \change{an agent is time-homogeneous, meaning that, in some sense, for all $i\in[N]$, the agents $(i,t)_{t\in[T]}$ are the same, and aggregates the stage costs across time.}
Consequently, the global aggregation operator $\operatorAggregation{\setProsumer \times \horizon}$ can be modeled as
aggregating, over \change{agents, their aggregated stage-costs}, \ie  $\operatorAggregation{\setProsumer \times \horizon} = \operatorAggregation{\setProsumer} \odot  \operatorTime{\horizon}{ }$ where the $\odot$ notation stands for
\begin{align}
    \operatorAggregation{\setProsumer} \odot  \operatorTime{\horizon}{ } \Bp{(\instantaneousCostProsumer{\indexProsumer}{\indexTime})_{\indexProsumer\in[\setProsumer],\indexTime\in[\horizon]}} 
    = \operatorAggregation{\setProsumer} \bgp{\operatorTime{\horizon}{1}\Bp{(\instantaneousCostProsumer{\indexTime}{1})_{\indexTime\in[\horizon]}}, \dots,\operatorTime{\horizon}{\setProsumer}\Bp{(\instantaneousCostProsumer{\indexTime}{\setProsumer})_{\indexTime\in[\horizon]}}}.
    \label{eq:composition}
\end{align}

Finally, we obtain the following model for the aggregation of \change{agents} in a dynamic framework:
\begin{subequations}
    \begin{align}
        \modelAggregatorDynamic{\horizon} \; := \quad \Min_{\decisionProsumerDynamic{\indexTime}{\indexProsumer}} \quad & 
        \operatorAggregation{\setProsumer} \odot  \operatorTime{\horizon}{ } \Bp{(\instantaneousCostProsumer{\indexTime}{\indexProsumer})_{\indexProsumer\in[\setProsumer],\indexTime\in[\horizon]}} \label{dynamic-obj2} \\
        \text{s.t.} \quad & \decisionProsumerDynamic{\indexTime}{\indexProsumer} \in \constraintProsumerDynamic{\indexTime}{\indexProsumer} & \forall \indexTime \label{feasible-space} \\
        & \change{ h(\decisionProsumerDynamic{\indexTime}{1}, \dots, \decisionProsumerDynamic{\indexTime}{\setProsumer})} \in \constraintMarketDynamic{\indexTime} & \forall \indexTime  \\
        & (\instantaneousCostProsumer{\indexTime}{\indexProsumer})_{\indexTime \in [\horizon]} \in \acceptableSet{ }{\indexProsumer}\label{eq:dynamic_acceptability},
    \end{align}
    \label{pb:aggregator_dynamic}
\end{subequations}
where we recall that we defined $\operatorAggregation{\horizon}$ as the sum, and suggest to 
choose $\operatorAggregation{\setProsumer}$ \change{from the fairness operators ($\objectiveUtilitarian{\setProsumer},\objectiveProportional{\setProsumer},\objectiveMinimaxSavings{\setProsumer}$) introduced in \cref{subsec:fairness_deterministic}.
Thus, we obtain a fair objective function of the aggregated model $\modelAggregatorDynamic{\horizon}$.
}
However, we have yet to adapt the notion of acceptability from \cref{subsec:acceptibility_deterministic} to this long-term framework, which is our next topic.


\subsection{Dynamic acceptability}
\label{subsec:dynamic_acceptability}

In long-term problems, for the aggregation to be acceptable, \change{agents} should not be tempted to leave the aggregation in between stages. 
Therefore, we extend our notion of acceptability constraint~\eqref{eq:acceptability} to a dynamic framework.
First, denote $\valueProsumerDynamic{\indexTime}{\indexProsumer}:=\instantaneousCostProsumer{\indexTime}{\indexProsumer}(\decisionProsumerDynamic{\indexTime}{\indexProsumer,\star})$, 
the optimal independent cost of \change{an agent} $\indexProsumer$ at stage $\indexTime$, 
where $\decisionProsumer{\indexProsumer,\star}$ is the optimal solution of 
Problem~\eqref{pb:prosumer_dynamic}.

The acceptability constraint \eqref{eq:acceptability} consist in requiring, for each agent $\indexProsumer$, that its vector of costs $(\instantaneousCostProsumer{\indexTime}{\indexProsumer})_{\indexTime\in[\horizon]}$ is less than $(\valueProsumerDynamic{\indexTime}{\indexProsumer})_{\indexTime\in[\horizon]}$.
Unfortunately, there is no natural ordering of $\RR^T$, and each (partial) order will define a different extension of the acceptability constraint~\eqref{eq:acceptability}.
We present now \change{some extensions of the acceptability constraint derived from standard partial orders}.

Maybe the most intuitive choice is the component-wise order (induced by the positive orthant), \ie comparing coordinate by coordinate.
This results in the \textit{stage-wise acceptability} constraint $\acceptableStagewise{ }$, 
which enforces that each agent benefits from the aggregation at each stage:
\begin{subequations}
\begin{align}
    \acceptableStagewise{\indexProsumer} = & \ba{ \; 
    (\utilityProsumerDynamic{t}{\indexProsumer})_{t\in[T]} \quad | \quad  
    \utilityProsumerDynamic{\indexTime}{\indexProsumer} \leq   \; \valueProsumerDynamic{\indexTime}{\indexProsumer}, \quad \forall \indexTime \in [\horizon]\; }. \label{acceptability:stagewise}
\end{align}
As this approach might be too conservative for our model, \new{\ie constrain the aggregation too much to take advantage of it}, we consider two other ordering. 

First, we can relax the stage-wise acceptability by considering that at each stage $\indexTime$, each agent benefits from the aggregation if we consider its costs aggregated up to time $\indexTime$.
This result in \textit{progressive acceptability} constraint $\acceptableProgressive{\indexProsumer}$:
\begin{align}
    \acceptableProgressive{\indexProsumer} = & \ba{\; (\utilityProsumerDynamic{t}{\indexProsumer})_{t\in[T]} 
    \quad | \quad
    \sum_{\tau=1}^\indexTime  \utilityProsumerDynamic{\tau}{\indexProsumer} \leq \; \sum_{\tau=1}^\indexTime \valueProsumerDynamic{\tau}{\indexProsumer}, \quad  \forall \indexTime \in [\horizon]\label{acceptability:progressive} \; } .
\end{align}

Second, we ensure that each agent, aggregating its cost over the whole horizon, benefits from the aggregation \change{(which amounts to the set in \eqref{eq:acceptability}) if we consider only the aggregated costs at the end of the horizon}.
We thus define the \textit{average acceptability} constraint:
\begin{align}
    \acceptableSet{av} {\indexProsumer} = & \ba{\; (\utilityProsumerDynamic{t}{\indexProsumer})_{t\in[T]} 
    \quad | \quad
    \sum_{\indexTime=1}^\horizon \utilityProsumerDynamic{\indexTime}{\indexProsumer} \leq \; \sum_{\indexTime=1}^\horizon \valueProsumerDynamic{\indexTime}{\indexProsumer} \; }.\label{acceptability:average} 
\end{align}
\end{subequations}

\begin{remark}
    We have that $\acceptableStagewise{\indexProsumer} \subseteq \acceptableProgressive{\indexProsumer} \subseteq \acceptableSet{av} {\indexProsumer}  $.
    The acceptability constraint should be chosen as to strike a balance between aggregated cost efficiency (obtained with \change{a less constrained} acceptability set), and incentive to stay in the aggregation (obtained with \change{a more constrained} acceptability set).
    \label{remark:bricks_dynamic_acceptability}
\end{remark}

\subsection{Numerical illustration}
\label{subsec:dynamic_illustration}

\begin{figure}[ht]
    \begin{adjustwidth}{-1.8cm}{}
    \centering
    \scalebox{1.1}{\pgfplotsset{compat=1.9}

\begin{tikzpicture}

\pgfplotsset{
        show sum on top/.style={
            /pgfplots/scatter/@post marker code/.append code={%
                \node[
                    at={(normalized axis cs:%
                            \pgfkeysvalueof{/data point/x},%
                            \pgfkeysvalueof{/data point/y})%
                    },
                    anchor=south,
                ]
                {\pgfmathprintnumber{\pgfkeysvalueof{/data point/y}}};
            },
            font=\footnotesize
        },
    }
    
 \begin{axis}[
     x=0.75cm,
     ybar stacked,
     bar width=0.5cm,
     width=0.8\linewidth,
     enlargelimits=true,
     ylabel={Total Costs (\$)},
     xtick={0, 1, 3, 4, 5, 6, 8, 9, 10, 11, 13, 14, 15, 16, 17},
     xticklabels={,\scriptsize $\modelProsumerDynamic{\horizon}{\indexProsumer}$, 
     \scriptsize $\modeName{U}{\emptyset}$,
     \scriptsize $\modeName{U}{av}$, 
     \scriptsize $\modeName{U}{p}$, 
     \scriptsize $\modeName{U}{s}$,
     \scriptsize $\modeName{SM}{\emptyset}$,
     \scriptsize $\modeName{SM}{av}$,
     \scriptsize $\modeName{SM}{p}$,
     \scriptsize $\modeName{SM}{s}$,
     \scriptsize $\modeName{P}{\emptyset}$,
     \scriptsize $\modeName{P}{av}$,  
     \scriptsize $\modeName{P}{p}$,  
     \scriptsize $\modeName{P}{s}$, },
      legend style={at={(0.5,0.95)},
        anchor=north,
        legend columns=4,
        draw=none
     },
     enlarge x limits={abs=-0.6cm}
 ]

\addplot +[fill=blue!40,draw=blue!70,nodes near coords, point meta=explicit symbolic, nodes near coords align={center}, nodes near coords style={blue, font=\scriptsize, anchor=center, yshift=0.1cm}] coordinates {(0,0) (1,50) (3,82) (4,26) (5,23) (6,10) (8,35) (9,35) (10,30) (11,10) (13,13) (14,13) (15,27) (16,10) (17,0)};

\addplot +[pattern = crosshatch, pattern color =red!50,draw=red!70, nodes near coords, point meta=explicit symbolic, nodes near coords align={center}, nodes near coords style={, red, font=\scriptsize, anchor=center, yshift=-0.1cm}] coordinates {(0,0) (1,280) (3,150) (4,175) (5,215) (6,240) (8,195) (9,195) (10,215) (11,240) (13,220) (14,220) (15,215) (16,240) (17,0)};

\addplot +[fill=brown!20,draw=brown!70, nodes near coords, point meta=explicit symbolic, nodes near coords align={center}, nodes near coords style={brown, font=\scriptsize, anchor=center, yshift=0.1cm,text opacity=1}] coordinates {(0,0) (1,40) (3,11) (4,40) (5,23) (6,8) (8,13) (9,13) (10,17) (11,8) (13,8) (14,8) (15,17) (16,8) (17,0)};

\addplot +[pattern = north west lines, pattern color=orange!70,draw=orange!70, nodes near coords, point meta=explicit symbolic, nodes near coords align={center}, nodes near coords style={orange, font=\scriptsize, anchor=center, yshift=-0.1cm}, show sum on top] coordinates {(0,0) (1,168) (3,90) (4,105) (5,129) (6,144) (8,117) (9,117) (10,129) (11,144) (13,131) (14,131) (15,129) (16,144) (17,0)};

\legend{$A_1$, $A_2$, $A_3$, $A_4$}
\end{axis}

\end{tikzpicture}}
    \end{adjustwidth}
    \caption{
    The \change{bars} correspond to the results of different models we solve.
    The first one is the independent model: we solve each $\modelProsumerDynamic{\horizon}{\indexProsumer}$ independently. 
    Then, there are three groups of four \change{bars}, each group corresponds to a choice of agent operator ($\objectiveUtilitarian{ },\objectiveMinimaxSavings{ },\objectiveProportional{ }$).
    Then, given an operator, we have the model first without then with different acceptability constraints ($\emptyset, \acceptableSet{av} { }, \acceptableProgressive{ }, \acceptableStagewise{ }$).  
    Each \change{bar} is decomposed in $4$ blocks corresponding to the share of each consumer $\indexProsumer$. \new{At the top of each bar, we can read the sum of aggregated costs in the corresponding model.}
    }
    \label{fig:dynamic_illustration}
\end{figure}

We take the same example as in \cref{subsec:deterministic_illustration} and try out different combinations of operator $\operatorAggregation{\setProsumer \times \horizon}$ ($\objectiveUtilitarian{ },\objectiveMinimaxSavings{ },\objectiveProportional{ }$) and acceptability set $\cA$ ($\emptyset, \acceptableSet{av} { }, \acceptableProgressive{ }, \acceptableStagewise{ }$).
\new{We denote $\modeName{f}{a}$ the model with agent operator $f \in \na{U,SM,P}$ and acceptability set $a \in \na{\emptyset, av, p, s}$.}
\Cref{fig:dynamic_illustration} represent the distribution of prosumers' costs for these different cases, while \cref{tab:results_dynamic} 
report their proportional savings. 
\new{Finally, on \cref{tab:details_purchases_dynamic}, we report the day-ahead and balancing purchases of the different models with progressive and stagewise acceptability, while the results with no acceptability and average acceptability are in \cref{subsec:deterministic_illustration} on \cref{tab:details_purchases}. }

\begin{table}[h]
\caption{Percentage of savings $\frac{\valueProsumer{\indexProsumer} - \objectiveProsumer{\indexProsumer}(\decisionProsumer{\indexProsumer})}{\valueProsumer{\indexProsumer}}$ achieved by $A_\indexProsumer$ in the model $\modeName{f}{a}$ depending on agent operator $f$ and acceptability set $a$ and $PoA$ (in percentage) of the corresponding model.}
\label{tab:results_dynamic}
\centering
\begin{tabular}{lccccc|ccccc|ccccc}
\toprule
& \multicolumn{5}{c}{Utilitarian $\objectiveUtilitarian{\setProsumer}$} & \multicolumn{5}{c}{Minimax $\objectiveMinimaxSavings{\setProsumer}$} & \multicolumn{5}{c}{Proportional $\objectiveProportional{\setProsumer}$} \\
\cmidrule(lr){2-6} \cmidrule(lr){7-11} \cmidrule(lr){12-16}
& A1 & A2 & A3 & A4  & $PoA$ & A1 & A2 & A3 & A4 & $PoA$ & A1 & A2 & A3 & A4  & $PoA$\\
\midrule
$\emptyset$ & -64  & 46  & 73  & 46 & 0 & 30  & 30  & 67  & 30 & 0 & 74  & 21  & 80  & 21 & 0  \\
$\acceptableSet{av} { }$ & 48  & 37  & 0  & 37 & 4 & 30  & 30  & 67  & 30 & 0 & 74  & 21  & 79  & 21 & 0 \\
$\acceptableProgressive{ }$ & 55 & 23 & 44 & 23 & 17 & 44 & 23 & 57 & 23 & 8.6 & 45 & 23 & 56 & 23 & 4.3 \\
$\acceptableStagewise{ }$ &  80 & 14 & 80 & 14 & 21 & 80 & 14 & 80 & 14 & 12 & 80 & 14 & 80 & 14 & 8.0 \\
\bottomrule
\end{tabular}
\end{table}

\new{
\begin{table}[ht]
\renewcommand{\arraystretch}{1.1}
\setlength{\tabcolsep}{4pt}

\begin{subtable}{.45\linewidth}
\centering
\begin{tabular}{l|cc|cc|cc|cc}
\hline
 & \multicolumn{2}{c}{A1} & \multicolumn{2}{c}{A2} & \multicolumn{2}{c}{A3} & \multicolumn{2}{c}{A4} \\ \hline
t & DA & B & DA & B &  DA & B & DA & B \\
1 & \; \textcolor{blue}{\textit{\textbf{0}}} \; & 0 & 0 & \textcolor{red}{\textit{\textbf{5}}} & \; \textcolor{blue}{\textit{\textbf{0}}} \; & 0 & 0 & \textcolor{red}{\textit{\textbf{3}}} \\
2 & \textcolor{blue}{\textit{\textbf{0}}} & 0 & 0 & \textcolor{red}{\textit{\textbf{5}}} & \textcolor{blue}{\textit{\textbf{0}}} & 0 & 0 & \textcolor{red}{\textit{\textbf{3}}} \\
3 & 5 & 0 & 5 & 0 & 4 & 0 & 3 & 0 \\
4 & 1.4 & 0 & 5 & 0 & 1.6 & 0 & 3 & 0 \\
5 & 3.6 & 0 & 5 & 0 & 2.4 & 0 & 3 & 0 \\
\hline
\end{tabular}
\caption{$\modeName{U}{\acceptableProgressive{ }}$}
\end{subtable}%
\begin{subtable}{.45\linewidth}
\centering
\begin{tabular}{||cc|cc|cc|cc}
    \hline
    
    \multicolumn{2}{c}{A1} & \multicolumn{2}{c}{A2} & \multicolumn{2}{c}{A3} & \multicolumn{2}{c}{A4} \\  \hline
    DA & B & DA & B &  DA & B &  DA & B \\
 \textcolor{blue}{\textit{\textbf{0}}} & 0 & 0 & \textcolor{red}{\textit{\textbf{5}}} & \textcolor{blue}{\textit{\textbf{0}}} & 0 & 0 & \textcolor{red}{\textit{\textbf{3}}} \\
 \textcolor{blue}{\textit{\textbf{0}}} & 0 & 0 & \textcolor{red}{\textit{\textbf{5}}} & 0 & \textcolor{blue}{\textit{\textbf{0}}} & 0 & \textcolor{red}{\textit{\textbf{3}}} \\
 5 & 0 & 5 & 0 & 4 & 0 & 3 & 0 \\
 \textcolor{blue}{\textit{\textbf{0}}} & 0 & 0 & \textcolor{red}{\textit{\textbf{5}}} & 0 & \textcolor{blue}{\textit{\textbf{0}}} & 0 & \textcolor{red}{\textit{\textbf{3}}} \\
 5 & 0 & 5 & 0 & 4 & 0 & 3 & 0 \\
\hline \end{tabular}
\caption{$\modeName{U}{\acceptableStagewise{ }}$ \textcolor{white}{$\modeName{U}{\acceptableProgressive{ }}$}}
\end{subtable}

\begin{subtable}{.45\linewidth}
\centering
    \begin{tabular}{l|cc|cc|cc|cc}
    \hline  & \multicolumn{2}{c}{A1} & \multicolumn{2}{c}{A2} & \multicolumn{2}{c}{A3} & \multicolumn{2}{c}{A4} \\  \hline
    t & DA & B & DA & B &  DA & B &  DA & B \\
1 & \; \textcolor{blue}{\textit{\textbf{0}}} \; & 0 & 0 & \textcolor{red}{\textit{\textbf{5}}} & \; \textcolor{blue}{\textit{\textbf{0}}} \; & 0 & 0 & \textcolor{red}{\textit{\textbf{3}}} \\
2 & \textcolor{blue}{\textit{\textbf{0}}} & 0 & 0 & \textcolor{red}{\textit{\textbf{5}}} & \textcolor{blue}{\textit{\textbf{0}}} & 0 & 0 & \textcolor{red}{\textit{\textbf{3}}} \\
3 & 5 & 0 & 5 & 0 & 4 & 0 & 3 & 0 \\
4 & 2 & 0 & 5 & 0 & 1 & 0 & 3 & 0 \\
5 & 3 & 0 & 5 & 0 & 3 & 0 & 3 & 0 \\
\hline \end{tabular}
\caption{$\modeName{SM}{\acceptableProgressive{ }}$}
\end{subtable}%
\begin{subtable}{.45\linewidth}
        \centering
        \begin{tabular}{||cc|cc|cc|cc}
        \hline  \multicolumn{2}{c}{A1} & \multicolumn{2}{c}{A2} & \multicolumn{2}{c}{A3} & \multicolumn{2}{c}{A4} \\  \hline
        DA & B & DA & B &  DA & B &  DA & B \\
 \textcolor{blue}{\textit{\textbf{0}}} & 0 & 0 & \textcolor{red}{\textit{\textbf{5}}} & \textcolor{blue}{\textit{\textbf{0}}} & 0 & 0 & \textcolor{red}{\textit{\textbf{3}}} \\
 \textcolor{blue}{\textit{\textbf{0}}} & 0 & 0 & \textcolor{red}{\textit{\textbf{5}}} & 0 & \textcolor{blue}{\textit{\textbf{0}}} & 0 & \textcolor{red}{\textit{\textbf{3}}} \\
 5 & 0 & 5 & 0 & 4 & 0 & 3 & 0 \\
 \textcolor{blue}{\textit{\textbf{0}}} & 0 & 0 & \textcolor{red}{\textit{\textbf{5}}} & 0 & \textcolor{blue}{\textit{\textbf{0}}} & 0 & \textcolor{red}{\textit{\textbf{3}}} \\
 5 & 0 & 5 & 0 & 4 & 0 & 3 & 0 \\
\hline \end{tabular}
\caption{$\modeName{SM}{\acceptableStagewise{ }}$ \textcolor{white}{$\modeName{U}{\acceptableProgressive{ }}$}}
\end{subtable}

\begin{subtable}{.45\linewidth}
    \centering
    \begin{tabular}{l|cc|cc|cc|cc}
    \hline  & \multicolumn{2}{c}{A1} & \multicolumn{2}{c}{A2} & \multicolumn{2}{c}{A3} & \multicolumn{2}{c}{A4} \\ \hline
    t & DA & B & DA & B &  DA & B &  DA & B \\
1 & \textcolor{blue}{\textit{\textbf{0}}} & 0 & 0 & \textcolor{red}{\textit{\textbf{5}}} & \textcolor{blue}{\textit{\textbf{0}}} & 0 & 0 & \textcolor{red}{\textit{\textbf{3}}} \\
2 & \textcolor{blue}{\textit{\textbf{0}}} & 0 & 0 & \textcolor{red}{\textit{\textbf{5}}} & \textcolor{blue}{\textit{\textbf{0}}} & 0 & 0 & \textcolor{red}{\textit{\textbf{3}}} \\
3 & 3.8 & 0 & 5 & 0 & 3.45 & 0 & 3 & 0 \\
4 & 1.94 & 0 & 5 & 0 & 1.06 & 0 & 3 & 0 \\
5 & 4.25 & 0 & 5 & 0 & 3.5 & 0 & 3 & 0 \\
\hline \end{tabular}
\caption{$\modeName{P}{\acceptableProgressive{ }}$}
\end{subtable}%
\begin{subtable}{.45\linewidth}
    \centering
    \begin{tabular}{||cc|cc|cc|cc}
    \hline  \multicolumn{2}{c}{A1} & \multicolumn{2}{c}{A2} & \multicolumn{2}{c}{A3} & \multicolumn{2}{c}{A4} \\ \hline
 DA & B & DA & B &  DA & B &  DA & B \\
 \textcolor{blue}{\textit{\textbf{0}}} & 0 & 0 & \textcolor{red}{\textit{\textbf{5}}} & \textcolor{blue}{\textit{\textbf{0}}} & 0 & 0 & \textcolor{red}{\textit{\textbf{3}}} \\
 \textcolor{blue}{\textit{\textbf{0}}} & 0 & 0 & \textcolor{red}{\textit{\textbf{5}}} & 0 & \textcolor{blue}{\textit{\textbf{0}}} & 0 & \textcolor{red}{\textit{\textbf{3}}} \\
 5 & 0 & 5 & 0 & 4 & 0 & 3 & 0 \\
 \textcolor{blue}{\textit{\textbf{0}}} & 0 & 0 & \textcolor{red}{\textit{\textbf{5}}} & 0 & \textcolor{blue}{\textit{\textbf{0}}} & 0 & \textcolor{red}{\textit{\textbf{3}}} \\
 5 & 0 & 5 & 0 & 4 & 0 & 3 & 0 \\
\hline \end{tabular}
\caption{$\modeName{P}{\acceptableStagewise{ }}$ \textcolor{white}{$\modeName{U}{\acceptableProgressive{ }}$}}
\end{subtable}
\caption{Quantity of energy purchased on the day-ahead and balancing markets per stage for all consumers depending on different models: in bold italic, we highlight purchases on the (more expensive) balancing market and stages where no purchases are made on the day-ahead market.}
\label{tab:details_purchases_dynamic}
\end{table}

}

\change{
We observe on  \cref{fig:dynamic_illustration} that increasing acceptability constraints (from none to average, progressive, and finally stage-wise) comes at a price but gives stronger guarantees to each prosumer. 
Indeed, we have seen that $\modeName{U}{\emptyset}$ is the most efficient model but yields solutions in contradiction with individual rationality. 
We can correct this defect by enforcing average acceptability, 
but this is not enough to ensure everyone gains from the aggregation, as $A_3$ makes $0\%$ of savings.
With more constrained acceptability, $\acceptableProgressive{ }$ and $\acceptableStagewise{ }$, we enforce individual rationality over time or at all times.
This leads to solutions where the savings among prosumers are shared in fairer proportions- at the loss of efficiency.
}

\change{
On the other hand, with an agent operator reflecting fairness (like scaled minimax or proportional), we obtain solutions that already aim at a fairer distribution of savings.
Consequently, if we can observe solutions changing with increasing acceptability constraints, those changes are clearer with a utilitarian operator.}
Indeed, in the utilitarian model, $A_2$ achieves savings ranging from $14\%$ to $46\%$ of his independent cost.
In contrast, under the scaled minimax approach, the savings range from $14\%$ to $30\%$,
and with the proportional approach, the savings fall between $14\%$ and $21\%$.

\change{
Note that even though the acceptability constraints and the agent operator are two distinct tools, they both drive the model to fairer solutions for all agents in the aggregation.
}

\section{Accommodating fairness to uncertainties with stochastic optimization}
\label{sec:stochastic}
Problems with energy generation, especially from renewable sources, and prices on energy markets are inherently uncertain \new{as we have decisions to make over time, and the future is uncertain}. 
Then, in addition to acceptability and fairness, we must tackle the challenge of handling uncertainties (while being fair about how we handle those). 
We want to address this issue by extending the problem presented in \cref{sec:problem} to a stochastic framework.
To that end, we introduce random variable $\va \xi$, along with probability space $\probabilitySpace$, which gathers all sources of uncertainties in the problem. 
\change{We assume that $\Omega$ is finite, a common simplification in stochastic optimization to make problems more tractable. If $\Omega$ is not finite we rely on sample average approximation.}

In the same way that we decomposed the problem in \cref{sec:dynamic} with $\nbStages$ time steps, 
we can decompose the problem here with $\Omega$ scenarios. 
Thus, there are similarities with the previous section. 
The main difference is that the set of time-step $\{1,\dots,T\}$ has a natural ordering,
while the set of scenario $\Omega$ does not, which leads to discussing different partial orders on $\RR^\Omega$ than on $\RR^T$.


\subsection{Static stochastic problem formulation}
\label{subsec:stochastic_problem}

The problem at hand is naturally formulated as a multi-stage stochastic problem.
For simplicity reasons, we first consider a $2-$stage relaxation of the problem:
in the first stage, \textit{here-and-now} decisions must be made before knowing the noise's realization;
in the second stage, once the noise's realization is revealed, \textit{recourse} actions can be decided.

We first adapt the individual model $\modelProsumer{\indexProsumer}$ to a stochastic framework:
\begin{subequations}
    \begin{align}
        \modelProsumerStochastic{ }{\indexProsumer}{\indexRiskMeasure} \; := \quad \min_{\decisionProsumerRandom{\indexProsumer}{\noiseVA}} \quad & \indexRiskMeasure \left[ \objectiveProsumer{\indexProsumer}(\decisionProsumerRandom{\indexProsumer}{\noiseVA}, \noiseVA) \right]  \\
        \text{s.t.} \quad & \decisionProsumerRandom{\indexProsumer}{\noiseVA} \in  \constraintProsumer{\indexProsumer} & \text{a.s.} \\
        & \change{\decisionProsumerRandom{\indexProsumer}{\noiseVA}} \in \constraintMarket & \text{a.s.},
    \end{align}
    \label{pb:prosumer_stochastic}
\end{subequations}
where $\indexRiskMeasure$ is a (coherent) risk-measure \ie a function which gives a deterministic cost equivalent to a random cost, reflecting the risk of a decision for prosumer $\indexProsumer$,
see \eg \cite{artzner1999coherent}.
The choice of $\indexRiskMeasure$ depends on the attitude of $\indexProsumer$ towards risk.
For example, the risk measure associated with a risk-neutral approach is the mathematical expectation $\EE_\xi$.
Alternatively, a highly risk-averse profile will opt for the worst-case measure $\sup_\xi$. 
Another widely used risk measure is the Average Value at Risk (a.k.a Conditional Value at Risk, or expected shortfall, see \cite{rockafellar2000optimization}), or a convex combination of expectation and Average Value at Risk.

Now, we adapt the deterministic aggregation model $\modelAggregator$.
We face the same challenge as in \cref{sec:dynamic}. 
With multiple scenarios, we can consider that we have $\setProsumer\times\Omega$ prosumers and we need to choose an operator  $\operatorAggregation{\setProsumer \times \Omega}: \RR^{\setProsumer \times \Omega} \to \RR$, leading to:
\begin{subequations}
    \begin{align}
        \modelAggregatorStochastic{ }{\indexRiskMeasure} \; := \quad \Min_{\va x} \quad & 
        \operatorAggregation{\setProsumer \times \Omega}
        \Bp{(\objectiveProsumer{\indexProsumer}(\decisionProsumerRandom{\indexProsumer}{\noiseVA}, \noiseVA))_{\indexProsumer \in [\setProsumer]}}  \\
        \text{s.t.} \quad & \decisionProsumerRandom{\indexProsumer}{\noiseVA} \in  \constraintProsumer{\indexProsumer} & \forall \indexProsumer  \in [\setProsumer]\quad \text{a.s.} \\
        & \change{h(\decisionProsumerRandom{1}{\noiseVA}, \dots, \decisionProsumerRandom{\setProsumer}{\noiseVA})} \in \constraintMarket & \text{a.s.} \\
        & \objectiveProsumer{\indexProsumer}(\decisionProsumerRandom{\indexProsumer}{\noiseVA}, \noiseVA) \in \acceptableSet{ }{\indexProsumer} & \forall \indexProsumer \in [\setProsumer] \quad \text{a.s.}. 
    \end{align}
    \label{pb:aggregation_simple_stochastic}
\end{subequations}
We \change{ assume the aggregator} knows risk measures and prosumers objectives. 
As in \cref{subsec:dynamic_problem}, there are multiple possible choices for such operators. 
We \delete{can} assume that this operator $\operatorAggregation{\setProsumer \times \Omega}$ results from the composition of two operators:
an uncertainty-operator $\riskMeasure{\Omega}{\indexProsumer}$ dealing with the scenarios, which can differ from one prosumer to another;
and an \operatorName $\operatorAggregation{\setProsumer}$, as defined in \cref{subsec:fairness_deterministic}.
However, contrary to \cref{sec:dynamic}, it is not clear if we should aggregate first with respect to uncertainty (meaning that a prosumer manages its own risk) or with respect to prosumers (meaning that the risks are shared).
We next discuss reasonable modeling choices of aggregation operators and acceptability constraints. 

\subsection{Stochastic objective}
\label{subsec:stochastic_objective}

For the sake of conciseness, we are going to consider two possible uncertainty aggregators: a risk-neutral choice, where $\riskMeasure{\Omega}{\indexProsumer}$ is the mathematical expectation $\riskNeutral{\noise}$,
and a worst-case operator where $\riskMeasure{\Omega}{\indexProsumer}$ is the supremum over the possible realization $\riskRobust{\noise}$.
For the \operatorName $\operatorAggregation{\setProsumer}$, which reflects the way to handle fairness, we consider either the utilitarian $\objectiveUtilitarian{\setProsumer}$ or the scaled minimax $\objectiveMinimaxSavings{\setProsumer}$ options (see \cref{subsec:acceptibility_deterministic} for definitions). 

We suggest four different compositions of $\riskMeasure{\Omega}{\indexProsumer}$ and $\operatorAggregation{\setProsumer}$ to construct the aggregation operator $\operatorAggregation{\setProsumer \times \Omega}$.
Again, for simplicity of notations, we write $\va{\objectiveProsumer{\indexProsumer}}$ instead of $\objectiveProsumer{\indexProsumer}(\decisionProsumerRandom{\indexProsumer}{\noiseVA}, \noiseVA)$.

First, we introduce the risk-neutral and utilitarian operator $\objectiveUtilitarianStochastic$, which aims at minimizing the aggregated expected costs of prosumers: 
\begin{subequations}
    \begin{align}
        \objectiveUtilitarianStochastic\Bp{(\va{\objectiveProsumer{\indexProsumer}})_{\indexProsumer \in [\setProsumer]}} & = \objectiveUtilitarian{\setProsumer} \odot \riskNeutral{\Omega} \Bp{(\va{\objectiveProsumer{\indexProsumer}})_{\indexProsumer \in [\setProsumer]}} \\
        & = \sum_{\indexProsumer=1}^\setProsumer \; \sum_{\noise \in \Omega} \; \pi_\noise \; \objectiveProsumer{\indexProsumer}(\decisionProsumerRealization{\indexProsumer}{\noise},\noise).
    \end{align}
\end{subequations}

Alternatively, considering a robust approach to uncertainties, we have the operator $\objectiveUtilitarianRobust$ which minimizes the worst-case aggregated costs of prosumers:
\begin{subequations}
    \begin{align}
        \objectiveUtilitarianRobust\Bp{(\va{\objectiveProsumer{\indexProsumer}})_{\indexProsumer \in [\setProsumer]}} & = \riskRobust{\noise \in \Omega} \odot \objectiveUtilitarian{\setProsumer} \Bp{(\va{\objectiveProsumer{\indexProsumer}})_{\indexProsumer \in [\setProsumer]}} \\
        & = \riskRobust{\noise \in \Omega} \bga{ \;\sum_{\indexProsumer=1}^\setProsumer \; \objectiveProsumer{\indexProsumer}(\decisionProsumerRealization{\indexProsumer}{\noise}, \noise) \;}.
    \end{align}
\end{subequations}

\vspace{2mm}
\begin{remark}
    We claim that $\riskRobust{\noise \in \Omega} \odot \objectiveUtilitarian{\setProsumer}$ makes more sense than $\objectiveUtilitarian{\setProsumer} \odot \riskRobust{\noise \in \Omega} $ as the later aggregates \change{each prosumer's worst costs.
    Indeed, if the worst-case costs for different prosumers occur in different scenarios, the aggregated costs calculated might never happen or happen in a scenario $\xi$ not in $\Omega$. 
    }

    On the other hand, we have \change{$\riskNeutral{\noise \in \Omega} \odot \objectiveUtilitarian{\setProsumer} = \objectiveUtilitarian{\setProsumer}\odot \riskNeutral{\noise \in \Omega}$}, by associativity of sums.
    Similarly, by associativity of supremum, we have $\riskRobust{\noise \in \Omega} \odot \objectiveMinimaxSavings{\setProsumer} = \objectiveMinimaxSavings{\setProsumer}
    \odot \riskRobust{\noise \in \Omega}$.
\end{remark}
\vspace{2mm}

As the first two operators do not model fairness considerations into the model,
we now look for a fair distribution by using $\objectiveMinimaxSavings{\setProsumer}$ to aggregate prosumers' costs. 
First, let $\decisionProsumerRandom{\indexProsumer,\star}{\noiseVA}$ be the\footnote{We assume uniqueness of a way of selecting an optimal solution, as in \cref{sec:problem})} optimal solution of $\modelProsumerStochastic{ }{\indexProsumer}{\indexRiskMeasure}$, and denote $\valueProsumerStochasticRealization{\indexProsumer}{\indexRiskMeasure}{\noise}~:=~\objectiveProsumer{\indexProsumer}(\decisionProsumerRealization{\indexProsumer,\star}{\noise}, \noise)$, the cost incurred by $\indexProsumer$ when operating alone under uncertainty realization $\noise$. 
Finally, $\valueProsumerStochasticRandom{\indexProsumer}{\indexRiskMeasure}$ is the random variable taking values $\valueProsumerStochasticRealization{\indexProsumer}{\indexRiskMeasure}{\noise}$ for the respective realization $\noise$.

\change{Results given in \cref{subsec:deterministic_illustration,subsec:dynamic_illustration} suggest} that the scaled minimax approach suits our problem more than the proportional approach.
Thus, in a stochastic framework, we propose the operator $\objectiveMinimaxStochastic$
:
\begin{subequations}
    \begin{align}
        \objectiveMinimaxStochastic\Bp{(\va{\objectiveProsumer{\indexProsumer}})_{\indexProsumer \in [\setProsumer]}} & = \objectiveMinimaxSavings{\setProsumer} \odot \riskNeutral{\Omega} \Bp{(\va{\objectiveProsumer{\indexProsumer}})_{\indexProsumer \in [\setProsumer]}} \\
        & = \max_{\indexProsumer \in [\setProsumer]} \bga{\; \frac{\EE \bc{\valueProsumerStochasticRandom{\indexProsumer}{\riskNeutral{ }}} - \sum_{\noise \in \Omega} \; \pi_\noise \;  \objectiveProsumer{\indexProsumer}(\decisionProsumerRealization{\indexProsumer}{\noise}, \noise)}{ \EE \bc{\valueProsumerStochasticRandom{\indexProsumer}{\riskNeutral{ }}}} \;  }.
    \end{align}
\end{subequations}
Finally, combining the robust and the scaled minimax approaches, we obtain the operator $\objectiveMinimaxRobust$, which focuses on the prosumer having the worst worst-case proportional costs:
\begin{subequations}
    \begin{align}
        \objectiveMinimaxRobust\Bp{(\va{\objectiveProsumer{\indexProsumer}})_{\indexProsumer \in [\setProsumer]}} & = \riskRobust{\noise \in \Omega} \odot \objectiveMinimaxSavings{\setProsumer} \Bp{(\va{\objectiveProsumer{\indexProsumer}})_{\indexProsumer \in [\setProsumer]}} \\
        & = \riskRobust{\noise \in \Omega} \bga{ \; \max_{\indexProsumer \in [\setProsumer]} \; \ba{ \; \frac{\valueProsumerStochasticRealization{\indexProsumer}{\riskNeutral{ }}{\noise} - \objectiveProsumer{\indexProsumer}(\decisionProsumerRealization{\indexProsumer}{\noise},\noise) }{\valueProsumerStochasticRealization{\indexProsumer}{\riskNeutral{ }}{\noise}} \; } \;}.
    \end{align}
\end{subequations}

\vspace{2mm}
\begin{remark}
    Note that here, depending on the \change{sense of the combination between uncertainty-operator and agent-operator,} we could have a model with different risk-measure profiles for the prosumers.  

    \delete{Further, as already pointed out, other (coherent) risk measures uncertainty-aggregator could be used. 
    Similarly, other agent-aggregators, as those presented in 
    might be relevant as well.
    }
\end{remark}
\vspace{2mm}

We now turn to extending the acceptability constraint~\eqref{eq:acceptability} to a stochastic setting. 

\subsection{Stochastic dominance constraints}
\label{subsec:stochastic_acceptability}

As in \cref{subsec:dynamic_acceptability}, to induce acceptability, we require that, for each prosumer $\indexProsumer$, its random cost $\objectiveProsumer{\indexProsumer}(\decisionProsumerRandom{\indexProsumer}{\noiseVA}, \noiseVA)$ is less than the random cost of the independent model $\valueProsumerStochasticRandom{\indexProsumer}{\riskNeutral{ }}
$.
Unfortunately, there is no natural ordering of random variable (or equivalently of $\RR^\Omega$), and each (partial) order will define a different extension of the acceptability constraint~\eqref{eq:acceptability}.

We now present four acceptability constraints, using various ordering on the space of random variable, leveraging the stochastic dominance theory (see \cite{DR2003} for an introduction in the context of stochastic optimization).
In this section, we give the mathematical expression of acceptability constraints, 
but mixed integer formulation can be found in \cref{ann:constraints}.

\begin{subequations}
In a very conservative perspective, we consider the almost-sure order, comparing random variables scenario by scenario: 
\begin{align}
    \acceptableAlmostSure{\indexProsumer} := \ba{\; \utilityProsumerStochasticRandom{\indexProsumer}{\indexRiskMeasure}  \, | \; \utilityProsumerStochasticRealization{\indexProsumer}{\indexRiskMeasure}{\indexScenario} \leq   \; \valueProsumerStochasticRealization{\indexProsumer}{\indexRiskMeasure}{\indexScenario}, & \quad \forall \indexScenario \; }.
    \label{eq:acceptability_almostsure}
\end{align}

We can relax the almost-sure ordering by not requiring the benefit of aggregation for all scenarios but distributionally. 
For example, if we have two scenarios $\xi$ and $\zeta$, with same probability, we consider that \change{it is acceptable to lose on $\xi$ if we do better on $\zeta$, that is such that $\utilityProsumerStochasticRealization{\indexProsumer}{\indexRiskMeasure}{\xi} \leq  \valueProsumerStochasticRealization{\indexProsumer}{\indexRiskMeasure}{\xi} $
and 
$ \utilityProsumerStochasticRealization{\indexProsumer}{\indexRiskMeasure}{\zeta} \geq  \valueProsumerStochasticRealization{\indexProsumer}{\indexRiskMeasure}{\zeta} $.}
To formalize this approach, we turn to \textit{stochastic first-order dominance constraints} (see \cite{DR2003}), and leverage $1^{st}$ order acceptability:
\begin{align}
   \acceptableFirstdc{\indexProsumer} := \quad & \ba{\; \utilityProsumerStochasticRandom{\indexProsumer}{\indexRiskMeasure}  \, | \; \utilityProsumerStochasticRandom{\indexProsumer}{\indexRiskMeasure}\preceq_{(1)}   \; \valueProsumerStochasticRandom{\indexProsumer}{\indexRiskMeasure} \; }\label{eq:acceptability_firstdc} \\
  := \quad & \ba{ \; \utilityProsumerStochasticRandom{\indexProsumer}{\indexRiskMeasure}  \, | \; \PP( \utilityProsumerStochasticRandom{\indexProsumer}{\indexRiskMeasure} > \eta) \; \leq \; \PP( \valueProsumerStochasticRandom{\indexProsumer}{\indexRiskMeasure} > \eta), \hspace{-1cm} & \forall \eta \in\RR \; }
  \notag\\
  := \quad & \ba{\; \utilityProsumerStochasticRandom{\indexProsumer}{\indexRiskMeasure}  \, |  \; \EE\bc{g(\utilityProsumerStochasticRandom{\indexProsumer}{\indexRiskMeasure})} 
  \; \leq \; 
  \EE \bc{g(\valueProsumerStochasticRandom{\indexProsumer}{\indexRiskMeasure})} \hspace{-1cm} & \forall g:\RR\to\RR, \text{ non-decreasing} }.
  \notag
\end{align}
One downside of this acceptability constraint is that the modeling entails numerous binary variables, posing practical implementation challenges.

We can thus consider a relaxed, less risk-averse version of $1^{st}$ order acceptability, relying on \textit{stochastic second-order dominance constraints}, also known as increasing convex acceptability, which is equivalent to :
\begin{align}
   \acceptableSeconddc{\indexProsumer} := \; & \ba{\; \utilityProsumerStochasticRandom{\indexProsumer}{\indexRiskMeasure}  \, | \;  \utilityProsumerStochasticRandom{\indexProsumer}{\indexRiskMeasure} \preceq_{(ic)} \; \valueProsumerStochasticRandom{\indexProsumer}{\indexRiskMeasure}\;  } \label{eq:acceptability_seconddc} \\
    = \; & \ba{\; \utilityProsumerStochasticRandom{\indexProsumer}{\indexRiskMeasure}  \, |  \; \EE\bc{\, ( \utilityProsumerStochasticRandom{\indexProsumer}{\indexRiskMeasure} - \eta)^+\, }\, \leq \EE\bc{\, ( \valueProsumerStochasticRandom{\indexProsumer}{\indexRiskMeasure} - \eta)^+ \,}
     & \forall \eta \in \RR}  \notag\\
   = \; & \ba{\; \utilityProsumerStochasticRandom{\indexProsumer}{\indexRiskMeasure}  \, | \; \EE\bc{\, g( \utilityProsumerStochasticRandom{\indexProsumer}{\indexRiskMeasure})\, }\, \leq \EE\bc{\, g( \valueProsumerStochasticRandom{\indexProsumer}{\indexRiskMeasure} ) \,}, & \hspace{-1.5cm}\forall g:\RR\to\RR, \text{convex, non-decreasing}\; }.
   \notag
\end{align} 
Moreover, increasing convex acceptability is also easier to implement than $1^{st}$ order acceptability (see \cref{ann:constraints}).

Finally, the risk-neutral acceptability constraint compares two random variables through their expectation:
\begin{align}
    \acceptableAverage{\indexProsumer} := \ba{\; \utilityProsumerStochasticRandom{\indexProsumer}{\indexRiskMeasure}  \, | \; \EE_\PP \, \nc{\utilityProsumerStochasticRandom{\indexProsumer}{\indexRiskMeasure}} \leq   \; \EE_\PP \, \nc{\valueProsumerStochasticRandom{\indexProsumer}{\indexRiskMeasure}} \;}. 
    \label{eq:acceptability_average}
\end{align}
\end{subequations}
We can use another convex risk measure instead of the expectation in \change{\eqref{eq:acceptability_average}}.

\vspace{2mm}
\begin{remark}
    We have that $\acceptableAlmostSure{\indexProsumer}~\subseteq~\acceptableFirstdc{\indexProsumer}~\subseteq~\acceptableSeconddc{\indexProsumer}~\subseteq~\acceptableAverage{\indexProsumer}$.
    \change{Therefore, in the same way as in \cref{remark:bricks_dynamic_acceptability}, the acceptability constraint yields a balance between risk-neutral ($ \acceptableAverage{ }$) and a robust approach on risk ($\acceptableAlmostSure{\indexProsumer}$), with intermediary visions on risk ($\acceptableSeconddc{\indexProsumer}$, $\acceptableFirstdc{\indexProsumer}$)}
    \label{remark:bricks_stochastic_acceptability}
\end{remark}
\vspace{2mm}




\subsection{Numerical illustration}
\label{subsec:stochastic_illustration}

We consider the stochastic version of the example presented in \cref{subsec:deterministic_illustration}, 
where balancing prices $\na{\pricesBRandomExample{\indexStage}}_{\indexStage \in [\horizon]}$ are random variables with uniform, independent, distribution over $\nc{0.35\pricesDAExample{\indexStage}, 5\pricesDAExample{\indexStage}}$.  
The problem can be formulated as a multi-stage program,
where day-ahead purchases are decided in the first stage \change{(a day-ahead of following stages)},
and then each stage corresponds to a time slot where we can buy energy on the balancing market at a price $\pricesBRandomExample{\indexStage}$. 

We solve and discuss the sample average approximation of the two-stage approximation of this problem. 
More precisely, we draw $50$ prices scenario, and solve a two-stage program where the first stage decisions are the day-ahead purchases, and the second stage decisions are the balancing purchases from time slot $1$ to $\horizon$.
We set $\setProsumer=4, \nbStages=10$, and we draw $\Omega=50$ scenarios of balancing prices.
For the prosumers' parameters and market prices, we use the data on \cref{tab:example_prices_stochastic,tab:parameters}.

\begin{table}[h]
    \caption{Prices on both markets}
\label{tab:example_prices_stochastic}
\centering
\def\arraystretch{1.5}
    \begin{tabular}{|c||*{10}{c|}}\hline
        $\indexStage$ & \textbf{1} & \textbf{2} & \textbf{3} & \textbf{4} & \textbf{5}& \textbf{6} & \textbf{7} & \textbf{8} & \textbf{9} & \textbf{10}\\ 
        \hline\hline
        $\pricesDAExample{\indexStage}$  & 3 & 3 & 7 & 4 & 2 & 10 & 7 & 4 & 7.5 & 8 \\\hline
        $\minDA{\indexStage}$ & 12 & 12 & 12 & 12 & 12 & 12 & 12 & 12 & 12 & 12 \\\hline
    \end{tabular}
\end{table}

\begin{table}[h]
\caption{Percentage of expected savings $\frac{\EE \bc{\valueProsumerStochasticRandom{\indexProsumer}{\riskRobust{ }}} - \EE \bc{ \objectiveProsumer{\indexProsumer}(\decisionProsumerRandom{\indexProsumer}{\noiseVA} }}{ \EE \bc{\valueProsumerStochasticRandom{\indexProsumer}{\riskRobust{ }}}}$ made by $A_\indexProsumer$, expected aggregated costs $\EE \bc{\operatorAggregation{\setProsumer}  (\objectiveProsumer{1}(\decisionProsumerRandom{\indexProsumer,\star}{\noiseVA})_{\indexProsumer \in [\setProsumer]})}$ and $PoA$ in the corresponding model.}
\label{tab:results}

\centering
\begin{tabular}{lcccccc|cccccc}
\toprule
& \multicolumn{6}{c}{Utilitarian Stochastic $\objectiveUtilitarianStochastic{}{}$} & \multicolumn{6}{c}{Scaled Minimax Stochastic $\objectiveMinimaxStochastic{}{}$} \\
\cmidrule(lr){2-7} \cmidrule(lr){8-13}
 & $A_1$ & $A_2$ & $A_3$ & $A_4$ & (A) & $PoA$ & $A_1$ & $A_2$ & $A_3$ & $A_4$ & (A) & $PoA$ \\
\midrule
$\emptyset$ & 2  & 52  & -14  & 52 & 684 & 0 & 32  & 36  & 32  & 32  & 770&  0\\
$\acceptableAverageExpected{ }$ & 4  & 51  & 1  & 48  & 684 & 0 & 32  & 36  & 32  & 32  & 770 &  0 \\
$ \acceptableSeconddcExpected{ }$ & 22  & 44  & 25  & 36  & 721 & 37 & 32  & 36  & 32  & 32  & 770& 0  \\
$\acceptableFirstdcExpected{ }$ & 36  & 24  & 34  & 18  & 882 & 198 & 23  & 21  & 28  & 19  & 918 & 148\\
$\acceptableAlmostSureExpected{ }$ & 43  & 17  & 41  & 14  & 930 & 246 & 37  & 17  & 33  & 16  & 941 & 171  \\
\toprule
& \multicolumn{6}{c}{Utilitarian Robust $\objectiveUtilitarianRobust{}{}$} & \multicolumn{6}{c}{Scaled Minimax Robust $\objectiveMinimaxRobust{}{}$} \\
\cmidrule(lr){2-7} \cmidrule(lr){8-13}
  & $A_1$ & $A_2$ & $A_3$ & $A_4$ & (A) & $PoA$ & $A_1$ & $A_2$ & $A_3$ & $A_4$ & (A) & $PoA$   \\
\midrule
$\emptyset$ & -19  & 53  & 6  & 49  & 686 & 0 & 22  & 48  & 21  & 38  & 693 & 0\\
$\acceptableAverageExpected{ }$ & 0  & 51  & 0  & 48  & 686 & 0 & 22  & 48  & 21  & 38  & 693 & 0\\
$ \acceptableSeconddcExpected{ }$ & 22  & 43  & 25  & 33  & 738 & 52 & 29  & 43  & 29  & 33  & 720 & 27 \\
$\acceptableFirstdcExpected{ }$ & 25  & 22  & 28  & 15   & 920 &234  & 28  & 26  & 33   & 18   & 881 & 188  \\
$\acceptableAlmostSureExpected{ }$ & 43  & 17  & 40  & 14  & 930 & 244& 43  & 17  & 40   & 14   & 930 & 237 \\
\bottomrule
\end{tabular}
\end{table}

We solve the problem with different combinations of aggregation operators and acceptability sets and can compare the impact of each combination on the solution. 
\new{We denote $\modeName{f}{a}$ the model with aggregation operator $f \in \na{US, SMS, UR, SMR}$ and acceptability set $a \in \na{\emptyset, \EE, (ic), (1), (a.s)}.$}
We read prosumers' expected percentage of savings with risk-neutral and worst-case approaches on \Cref{tab:results}. 
\new{For example, in model $\modeName{US}{(ic)}$, we read that $\consumer{1}$ (resp. $\consumer{2}, \consumer{3}, \consumer{4}$) saves $22\%$ (resp. $44\%, 25\%, 36\%$) of its costs.
The expected cost of the aggregation is $721 \$ $, thus asking for increasing-convex acceptability costs $37 \$$.
}
Moreover, we can observe the distribution of prosumers' expected costs with a risk-neutral (resp. worst-case) approach on \cref{fig:stochastic_illustration} (resp. \cref{fig:robust_illustration}).

Our first comment is that the problems previously identified from a utilitarian perspective with no acceptability constraints are still present in a stochastic framework. 
Indeed, both with the risk-neutral utilitarian $\objectiveUtilitarianStochastic$ and worst-case utilitarian $\objectiveUtilitarianRobust$ operators, 
we observe on \cref{tab:results} that some prosumers can pay more in the aggregation compared to being alone ($A_3$ pays $+14\%$ in the stochastic approach, and $A_1$ pays $+19\%$ in the robust approach).
This highlights the necessity for either acceptability constraints or an aggregation operator. 

\begin{figure}[ht]
    \begin{adjustwidth}{-2.1cm}{}
    \scalebox{1.1}{\pgfplotsset{compat=1.9}

\begin{tikzpicture}

\pgfplotsset{
        show sum on top/.style={
            /pgfplots/scatter/@post marker code/.append code={%
                \node[
                    at={(normalized axis cs:%
                            \pgfkeysvalueof{/data point/x},%
                            \pgfkeysvalueof{/data point/y})%
                    },
                    anchor=south,
                ]
                {\pgfmathprintnumber{\pgfkeysvalueof{/data point/y}}};
            },
            font=\footnotesize
        },
    }
    
 \begin{axis}[
     x=0.9cm,
     ybar stacked,
     bar width=0.5cm,
     width=0.8\linewidth,
     enlargelimits=true,
     ylabel={Total Expected Costs (\$)},
     xtick={1, 3, 4, 5, 6, 7, 9, 10, 11, 12, 13},
     xticklabels={$\modelProsumer{\indexProsumer}$, 
     $\modeName{US}{\emptyset}$,
     $\modeName{US}{\EE}$,
     $\modeName{US}{(ic)}$, 
     $\modeName{US}{(1)}$,  
     $\modeName{US}{a.s}$, 
     $\modeName{SMS}{\emptyset}$,
     $\modeName{SMS}{\EE}$,
     $\modeName{SMS}{(ic)}$, 
     $\modeName{SMS}{(1)}$,  
     $\modeName{SMS}{a.s}$,}, 
      legend style={at={(0.5,0.95)},
        anchor=north,
        legend columns=4,
        draw=none
     },
 ]

\addplot +[fill=blue!40,draw=blue!70,nodes near coords, point meta=explicit symbolic, nodes near coords align={center}, nodes near coords style={blue, font=\scriptsize, anchor=center, yshift=0.1cm}] coordinates { (1,90) (3,88) (4,86) (5,70) (6,57) (7,51) (9,61) (10,61) (11, 61) (12,69) (13,57) };

\addplot +[pattern = crosshatch, pattern color =red!50,draw=red!70, nodes near coords, point meta=explicit symbolic, nodes near coords align={center}, nodes near coords style={, red, font=\scriptsize, anchor=center, yshift=-0.1cm}] coordinates { (1,665) (3,316) (4,324) (5,371) (6,503) (7,553) (9,425) (10,425) (11, 425) (12,522) (13,553) };

\addplot +[fill=brown!20,draw=brown!70, nodes near coords, point meta=explicit symbolic, nodes near coords align={center}, nodes near coords style={brown, font=\scriptsize, anchor=center, yshift=0.1cm,text opacity=1}] coordinates { (1,120) (3,137) (4,118) (5,89) (6,78) (7,71) (9,81) (10,81) (11, 81) (12,86) (13,80) };

\addplot +[pattern = north west lines, pattern color=orange!70,draw=orange!70, nodes near coords, point meta=explicit symbolic, nodes near coords align={center}, nodes near coords style={orange, font=\scriptsize, anchor=center, yshift=-0.1cm}, show sum on top] coordinates { (1,297) (3,144) (4,156) (5,190) (6,243) (7,255) (9,202) (10,202) (11, 202) (12,241) (13,251) };

\legend{$A_1$, $A_2$, $A_3$, $A_4$}
\end{axis}

\end{tikzpicture}}
    \end{adjustwidth}
        \caption{
        The \change{bars} correspond to the results of different models we solve with a stochastic approach.
        The first one is the \change{model without aggregation}: we solve each $\modelProsumerStochastic{ }{\indexProsumer}{\riskNeutral{ }}$ independently. 
        The second \change{bar} corresponds to the problem solved with operator $\objectiveUtilitarianStochastic$ without acceptability constraints. 
        Then, the four following \change{bars} correspond to the same problem with increasingly strong acceptability ($ \acceptableAverageExpected{ }, \acceptableSeconddcExpected{ }, \acceptableFirstdcExpected{ },\acceptableAlmostSureExpected{ }$). 
        The following \change{bar} is for the problem solved with operator $\objectiveMinimaxStochastic$ without acceptability constraints, followed by four \change{bars} with different acceptability sets.
        Each \change{bar} is decomposed in $4$ blocks corresponding to the expected share $\EE \nc{\objectiveProsumer{\indexProsumer}(\decisionProsumerRandom{\indexProsumer}{\noiseVA}, \noiseVA)}$ of each consumer $\indexProsumer$.
        \new{At the top of each bar, we can read the sum of expected aggregated costs in the corresponding model.}
    }
    \label{fig:stochastic_illustration}
\end{figure}

If we choose a fair approach through the objective (operators $\objectiveMinimaxStochastic$ and $\objectiveMinimaxRobust$), 
we guarantee a higher percentage of savings to all prosumers than in the utilitarian approach. 
For example, with no acceptability constraints, all prosumers save at least $32\%$ of their costs in a risk-neutral approach and $21\%$ in a robust approach,
compared to respectively $-14\%$ and $-19\%$ with the utilitarian approach. 
This comes at the price of efficiency, especially in the risk-neutral case, as the expected aggregated costs of the scaled minimax approach is $13\%$ higher than with the utilitarian approach.
This remains true as \change{we use smaller acceptability set}. 

Conversely, when solving this problem with a utilitarian approach (operators $\objectiveUtilitarianStochastic$ and $\objectiveUtilitarianRobust$), we can increase the guaranteed percentage of savings by constraining more the acceptability.
Indeed, with $\objectiveUtilitarianStochastic$, all prosumers save at least from $1\%$ with expected acceptability to $22\%$ with increasing convex acceptability, 
and with $\objectiveUtilitarianRobust$, it is from $0\%$ to $22\%$.
However, increasing the acceptability to first-order or almost sure does not improve this guarantee, as the problem is now getting too constrained.
For example, with almost-sure acceptability, the choice of the operator on uncertainty is inconsequential:
the distribution of costs is the same with both operators $\objectiveUtilitarianStochastic$ and $\objectiveUtilitarianRobust$.
Notably, there exists a substantial gap between increasing-convex acceptability and first-order acceptability.
For example, with the scaled minimax stochastic operator $\objectiveMinimaxStochastic$, the costs increase from $770$ with increasing convex acceptability to $918$ with first-order acceptability.

\begin{figure}[ht]
    \begin{adjustwidth}{-2.2cm}{}
    \scalebox{1.1}{\pgfplotsset{compat=1.9}

\begin{tikzpicture}

\pgfplotsset{
        show sum on top/.style={
            /pgfplots/scatter/@post marker code/.append code={%
                \node[
                    at={(normalized axis cs:%
                            \pgfkeysvalueof{/data point/x},%
                            \pgfkeysvalueof{/data point/y})%
                    },
                    anchor=south,
                ]
                {\pgfmathprintnumber{\pgfkeysvalueof{/data point/y}}};
            },
            font=\footnotesize
        },
    }
    
 \begin{axis}[
     x=0.9cm,
     ybar stacked,
     bar width=0.5cm,
     width=0.9\linewidth,
     enlargelimits=true,
     ylabel={Total Expected Costs (\$)},
     xtick={1, 3, 4, 5, 6, 7, 9, 10, 11, 12, 13},
          xticklabels={$\modelProsumer{\indexProsumer}$, 
     $\modeName{UR}{\emptyset}$,
     $\modeName{UR}{\EE}$,
     $\modeName{UR}{(ic)}$, 
     $\modeName{UR}{(1)}$,  
     $\modeName{UR}{a.s}$, 
     $\modeName{SMR}{\emptyset}$,
     $\modeName{SMR}{\EE}$,
     $\modeName{SMR}{(ic)}$, 
     $\modeName{SMR}{(1)}$,  
     $\modeName{SMR}{a.s}$,}, 
      legend style={at={(0.5,0.95)},
        anchor=north,
        legend columns=4,
        draw=none
     },
 ]

\addplot +[fill=blue!30,draw=blue!40,nodes near coords, point meta=explicit symbolic, nodes near coords align={center}, nodes near coords style={blue, font=\scriptsize, anchor=center, yshift=0.1cm}] coordinates {  (1,90) (3,106) (4,90) (5,70) (6,67) (7,51) (9,70) (10,70) (11, 63) (12,65) (13,51)  };

\addplot +[pattern = crosshatch, pattern color =red!50,draw=red!70, nodes near coords, point meta=explicit symbolic, nodes near coords align={center}, nodes near coords style={, red, font=\scriptsize, anchor=center, yshift=-0.1cm}] coordinates {  (1,665) (3,315) (4,324) (5,381) (6,515) (7,553) (9,345) (10,345) (11, 376) (12,493) (13,553)  };

\addplot +[fill=brown!20,draw=brown!70, nodes near coords, point meta=explicit symbolic, nodes near coords align={center}, nodes near coords style={brown, font=\scriptsize, anchor=center, yshift=0.1cm,text opacity=1}] coordinates {  (1,120) (3,112) (4,120) (5,89) (6,86) (7,71) (9,94) (10,94) (11, 84) (12,80) (13,71)  };

\addplot +[pattern = north west lines, pattern color=orange!70,draw=orange!70, nodes near coords, point meta=explicit symbolic, nodes near coords align={center}, nodes near coords style={orange, font=\scriptsize, anchor=center, yshift=-0.1cm}, show sum on top] coordinates {  (1,297) (3,153) (4,153) (5,198) (6,258) (7,255) (9,183) (10,183) (11, 198) (12,244) (13,255)  };

\legend{$A_1$, $A_2$, $A_3$, $A_4$}
\end{axis}

\end{tikzpicture}}
    \end{adjustwidth}
        \caption{
        This figure can be read like \cref{fig:stochastic_illustration}, except that the two considered operators are $\objectiveUtilitarianRobust$ and $\objectiveMinimaxRobust$.
    }
    \label{fig:robust_illustration}
\end{figure}

Thus, we obtain various solutions with different balances between efficiency and fairness and different risk visions. 
\change{In this example, in the stochastic case, if we want to give the same guarantees to every prosumer, the natural choice is operator $\objectiveMinimaxStochastic$. 
However, this approach costs $13\%$ more than in $\modeName{\objectiveUtilitarianStochastic}{\emptyset}$. 
If we want to opt for a less costly approach, 
the operator $\objectiveUtilitarianStochastic$ combined with increasing convex acceptability might be considered a better trade-off between efficiency and fairness: 
it ensures at least $22\%$ of savings to each prosumer and induces $5\%$ of efficiency loss.  }
\new{Finally, increasing the variability of the scenarios drawn for these tests does not significantly impact the empirical conclusions we make in this section.}

\section*{Conclusion}
\label{sec:conclusion}
\new{   
In this paper, we have provided a framework with tools to accommodate fairness in prosumer aggregation problems.
Through the discussion in Section 2, we emphasized the importance of fairness and the need to carefully consider how to model it and be aware of the different approaches available.
Since modeling fairness aims to improve the social acceptability of mathematical models, we connected underlying philosophical concepts to the fairness modeling process.
}

\delete{We have introduced, developed, and analyzed models for fair prosumer aggregation.}
First, we discussed acceptability constraints to discourage prosumers from leaving the aggregation. We then compared different choices of the objective function (utilitarian, scaled minimax, and proportional).
Through a stylized (and \change{more straightforward to interpret}) deterministic case study, we showed how different combinations of objectives and constraints influence solutions, emphasizing the importance of fairness and acceptability considerations.
We then extended the model to dynamic and stochastic frameworks, aligning it with what we expect practical problems to be \delete{ (\ie, decision-making most likely is for a sequence of time periods, under uncertainty)}. 
In this context, we adapt acceptability constraints to account for long-term horizons and uncertainties, 
and we showcase their impact on solutions using similar stylized instances.

In our numerical example,
\new{
we obtained a spectrum of options from various combinations of acceptability sets and objective functions, ranging from the most efficient models (with the lowest aggregated costs) to the fairest models (where agents' gains are more comparable). 
Too-restrictive acceptability sets or a bargaining approach (proportional operator $\objectiveProportional{\setProsumer}$) can significantly reduce efficiency, while an intermediate approach leverages aggregation benefits without excessively favoring certain prosumers.}
\change{Thus, we recommend the proportional min-max agent aggregator with progressive acceptability constraint in the dynamic case (resp.\ increasing convex acceptability in the stochastic case), which balances efficiency and fairness well.}
Recall that the framework discussed here is not reduced to prosumer aggregation in energy markets only, and can be adapted to other aggregation problems in energy system management problems (\eg virtual power plant, portfolio management in energy markets, ancillary service provision, etc.).

In future work, we plan to discuss the extension of the aggregation problem to a multistage stochastic program, 
\change{where we would have to combine the partial orders presented in the dynamic framework in \cref{sec:dynamic} and the stochastic orders of the stochastic framework in \cref{sec:stochastic}}.
This will require a discussion of possible aggregators $\operatorAggregation{\horizon\times \Omega \times \setProsumer}$ over agent, time, and uncertainty simultaneously.
Although we can easily assume a factorization of the form $\operatorAggregation{\setProsumer} \odot \operatorAggregation{\horizon\times \Omega}$, it would not be realistic to describe $\operatorAggregation{\horizon\times \Omega}$ as
the composition of a time aggregator and an uncertainty aggregator.
Indeed, such a factorization would not guarantee time-consistency of the problem and might not even preserve non-anticipativity. 
Further, acceptability constraints must be defined using multivariate stochastic order (see \cite{dentcheva2009optimization,armbruster2015models,dentcheva2016two}) whose mathematical programming representations are more involved.
\new{Finally, it would be interesting to investigate potential bounds on \textit{price of acceptability}.}



\subsection*{Statements and Declarations}
The authors declare they have no financial interests.

\bibliography{bibliography}

\newpage

\begin{appendix}

\section{Computing Shapley's values}
\label{ann:shapley}

In this appendix, we give more details on how to compute Shapley's values for the application in \cref{subsec:deterministic_illustration}. 
First, we introduce the general formulas and definitions required to compute Shapley's values, 
then we apply them to our example.

\subsection{Definitions and Formulas}

Shapley's values are commonly seen as a fair distribution of costs (or revenues) when agents cooperate in a cooperative game. 
The Shapley Value returns each player's fair share of the total gains by averaging their contributions across all possible ways they can join the coalition.
Let $N$ be the number of agents in the game.
We denote $w: 2^{|N|} \rightarrow \RR$, the worth function associating to a coalition $S$, the expected payoff obtained by cooperation.

To compute Shapley's values, first, we compute $\delta_i(S)$, the marginal contribution of agent $i$ to coalition $S \subset N$:
\begin{subequations}
 \begin{align}
     \delta_i(S) = w(S \cup \na{i}) - w(S)
 \end{align}   
Then, the Shapley value of agent $i$, given a characteristic function $w$, is $\phi_i(w)$:

 \begin{align}
     \phi_i(w) = \frac{1}{n} \sum_{S \subset N \setminus \na{i}} \begin{pmatrix}
  n-1 \\ 
  |S| 
\end{pmatrix}^{-1} \delta_i(S) 
 \end{align}  
\end{subequations}

\subsection{Application to our example}

We consider an example with $4$ agents.
For each coalition $S \subset [4]$, we solve the aggregated problem with utilitarian operator $\objectiveUtilitarian{S}$ and acceptability constraints $\acceptableSet{1}{ }$-as in cooperation games we must satisfy individual rationality properties- and obtain optimal solution $c_S$
Then, we introduce the characteristic function assessing the worth of coalition $S$ as: 
\begin{subequations}
 \begin{align}
     w(S) = \sum_{i \in S} \valueProsumer{i} - c_S,
 \end{align}   
where $\valueProsumer{i}$ is the optimal value of problem $\modelProsumer{i}$.
Thus $w(S)$ designates the savings made by coalition $S$ when they cooperate.
On \cref{tab:intermediate,tab:intermediate2}, we detail the intermediate computations, and in \cref{eq:shapley_value1,eq:shapley_value2,eq:shapley_value3,eq:shapley_value4} we compute shapley's values. 

\begin{table}[ht]
\setlength{\tabcolsep}{2pt}
\renewcommand{\arraystretch}{1.2}
\centering
\begin{tabular}{c|cccc|cccccc|cccc|c|}
& \na{1} & \na{2} & \na{3} & \na{4} & \na{1,2} & \na{1,3} & \na{1,4} & \na{2,3} & \na{2,4} & \na{3,4} & \na{1,2,3} & \na{1,2,4} & \na{1,3,4} & \na{2,3,4} & \na{1,2,3,4}  \\\hline
$\underset{i \in S}{\sum} \valueProsumer{i}$ & 50  & 280  & 40  & 168  & 330  & 90  & 218  & 320  & 448  & 208  & 370  & 498  & 258  & 488  & 538 \\\hline 
$c_S$ & 50  & 280  & 40  & 168  & 330  & 90  & 218  & 320  & 448  & 208  & 244  & 365  & 162  & 392  & 333  \\ \hline 
$w(S)$ & 0  & 0  & 0  & 0  & 0  & 0  & 0  & 0  & 0  & 0  & 126  & 133  & 96  & 96  & 205  \\
\end{tabular}
\caption{Costs with and without cooperation and worth $w$ of each coalition $S$.}
\label{tab:intermediate}
\end{table}

\begin{table}[ht]
\setlength{\tabcolsep}{2pt}
\renewcommand{\arraystretch}{1.2}
\centering
\begin{tabular}{c|cccc|cccccc|cccc|}
& \na{1} & \na{2} & \na{3} & \na{4} & \na{1,2} & \na{1,3} & \na{1,4} & \na{2,3} & \na{2,4} & \na{3,4} & \na{1,2,3} & \na{1,2,4} & \na{1,3,4} & \na{2,3,4} \\\hline
$\delta_1(S)$ & -  & 0  & 0  & 0 & - & - & - & 126 & 133 & 96 & - & - & - & 109 \\\hline 
$\delta_2(S)$ & 0  & -  & 0  & 0  & - & 126 & 133  & -  & -  & 96 & - & - & 109  & - \\ \hline 
$\delta_3(S) $ & 0  & 0  & -  & 0  & 126  & -  & 96 & - & 96 & -  & - & 72  & -  & - \\\hline
$\delta_4(S)$ & 0  & 0  & 0  & -  & 133  & 96  & -  & 96  & -  & -  & 79  & -  & -  & - \\
\end{tabular}
\caption{Marginal contributions of agent $i$ to each coalition.}
\label{tab:intermediate2}
\end{table}

\begin{align}
    \phi_1(w) & = \frac{1}{4} \bgc{  
    \begin{pmatrix} 3 \\ 1 \end{pmatrix}^{-1} \bp{\delta_1(\na{2})+\delta_1(\na{3})+\delta_1(\na{4})} \notag \\
    & \qquad \qquad +
    \begin{pmatrix} 3 \\ 2 \end{pmatrix}^{-1} \bp{\delta_1(\na{2,3})+\delta_1(\na{2,4})+\delta_1(\na{3,4})}\notag \\
    & \qquad \qquad +
     \begin{pmatrix} 3 \\ 3 \end{pmatrix}^{-1} \delta_1(\na{2,3,4})} \notag \\
    \phi_1(w) & = \frac{1}{4} \bgc{\frac{1}{3} \times 0 + \frac{1}{3} (126+133+96) + 1 \times 109} =  56.8 \label{eq:shapley_value1} \\
   \phi_2(w) & = \frac{1}{4} \bgc{\frac{1}{3} \times 0 + \frac{1}{3} (126+133+96) + 1 \times 109} =  56.8 \label{eq:shapley_value2}\\ 
   \phi_3(w) & = \frac{1}{4} \bgc{\frac{1}{3} \times 0 + \frac{1}{3} (126+96 +96) + 1 \times 72} = 44.5 \label{eq:shapley_value3}\\
    \phi_4(w) & =  \frac{1}{4} \bgc{\frac{1}{3} \times 0 + \frac{1}{3} (133+96 +96) + 1 \times 79} = 46.8 \label{eq:shapley_value4}
\end{align}

\end{subequations}

As the values computed here represent the way to distribute savings, we must subtract $\phi_i(w)$ from the cost of agent $i$ to obtain its costs in the aggregation after fair allocation through Shapley's values in \cref{tab:shapley}.

\begin{table}[h]
\setlength{\tabcolsep}{6pt}
\renewcommand{\arraystretch}{1.5}
    \centering
    \begin{tabular}{c|c|c|c|c|}
         & $A_1$ & $A_2$ & $A_3$ & $A_4$ \\ \hline
        $\phi_i(w)$ &  56.8 & 56.8 & 44.5 & 46.8 \\ \hline 
        $v^i$ & 50 & 280 & 40 & 168 \\ \hline 
        $L^i(x^i)$ & -6.8 & 223.2 & -4.5 &  121.2 \\ \hline 
        $\frac{100(\valueProsumer{\indexProsumer} - \objectiveProsumer{\indexProsumer}(\decisionProsumer{\indexProsumer}))}{\valueProsumer{\indexProsumer}}$ & 114 & 20 & 111 & 28
    \end{tabular}
    \caption{Costs and proportional savings of the agents with Shapley's post-allocation scheme.}
    \label{tab:shapley}
\end{table}


\section{Modeling of stochastic dominance constraints}
\label{ann:constraints}

We present here practical formulas to implement the stochastic orders dominance constraints introduced in \cref{subsec:stochastic_acceptability}.
Those constraints establish a dominance between $\valueProsumerStochasticRandom{\indexProsumer}{\indexRiskMeasure}$, the random variable representing $\indexProsumer$ independent costs,
and $\utilityProsumerStochasticRandom{\indexProsumer}{\indexRiskMeasure}$, the random variable representing $\indexProsumer$ costs in the aggregation. 

\subsection{First-order dominance constraint model}
\label{ann:1_dominance_constraints}

The first-order dominance constraints \eqref{eq:acceptability_firstdc} model is based on \cite{GNS2008}.
\begin{lem}
    In Problem $\modelAggregatorStochastic{ }{\indexRiskMeasure}$, acceptability constraints $\utilityProsumerStochasticRandom{\indexProsumer}{\indexRiskMeasure}\preceq_{(1)} \; \valueProsumerStochasticRandom{\indexProsumer}{\indexRiskMeasure}$  can be modeled with: 
    \begin{subequations}
    \begin{align}
        & b^\indexProsumer_{\indexScenario,\indexScenarioDominance} \in \na{0,1} & \forall \indexScenarioDominance \in [ \nbScenario], \; \forall \indexScenario \in \Omega
        \\
        & \utilityProsumerStochasticRealization{\indexProsumer}{\indexRiskMeasure}{\indexScenario} - \valueProsumerStochasticRealization{\indexProsumer}{\indexRiskMeasure}{\indexScenarioDominance} \leq M \; b^\indexProsumer_{\indexScenario,\indexScenarioDominance} & \forall \indexScenarioDominance \in [ \nbScenario], \; \forall \indexScenario \in \Omega \\
        & \underset{\indexScenario=1}{\overset{\nbScenario}{\sum}} \pi_{\indexScenario} b^\indexProsumer_{\indexScenario,\indexScenarioDominance} \leq a_{\indexScenarioDominance} & \forall \indexScenarioDominance \in [ \nbScenario].
    \end{align}
    \end{subequations}
    We denote $a_{\indexScenarioDominance} :=\PP( \valueProsumerStochasticRandom{\indexProsumer}{\indexRiskMeasure} > \valueProsumerStochasticRealization{\indexProsumer}{\indexRiskMeasure}{\indexScenarioDominance})$, which is a parameter for the aggregation problem. 
\end{lem}

\begin{proof}
    
As $\Omega$ is assumed to be finite, $\valueProsumerStochasticRandom{\indexProsumer}{\indexRiskMeasure}$ follows discrete distribution with realizations $\valueProsumerStochasticRealization{\indexProsumer}{\indexRiskMeasure}{\indexScenarioDominance}$ for $\indexScenarioDominance \in \Omega$. 
Then,
\begin{align*}
    \utilityProsumerStochasticRandom{\indexProsumer}{\indexRiskMeasure}\preceq_{(1)} \; \valueProsumerStochasticRandom{\indexProsumer}{\indexRiskMeasure}
    \iff \quad & \PP( \utilityProsumerStochasticRandom{\indexProsumer}{\indexRiskMeasure} > \eta) \; \leq \; \PP( \valueProsumerStochasticRandom{\indexProsumer}{\indexRiskMeasure} > \eta) & \forall \eta \in \RR \\
    \iff \quad & \PP( \utilityProsumerStochasticRandom{\indexProsumer}{\indexRiskMeasure} > \valueProsumerStochasticRealization{\indexProsumer}{\indexRiskMeasure}{\indexScenarioDominance}) \; \leq \; \PP( \valueProsumerStochasticRandom{\indexProsumer}{\indexRiskMeasure} > \valueProsumerStochasticRealization{\indexProsumer}{\indexRiskMeasure}{\indexScenarioDominance}) & \forall \indexScenarioDominance \in \Omega. 
\end{align*}
Then, 
using $\PP(\va X > x) = \EE[\1_{\va X > x}]$, and
 introducing binary variables $b^\indexProsumer_{\indexScenario,\indexScenarioDominance}=\1_{ \utilityProsumerStochasticRealization{\indexProsumer}{\indexRiskMeasure}{\indexScenario} > \valueProsumerStochasticRealization{\indexProsumer}{\indexRiskMeasure}{\indexScenarioDominance}}$, 
we get:
\begin{align*}
    \bgp{\PP( \utilityProsumerStochasticRandom{\indexProsumer}{\indexRiskMeasure} > \valueProsumerStochasticRealization{\indexProsumer}{\indexRiskMeasure}{\indexScenarioDominance}) \; \leq \; \PP( \valueProsumerStochasticRandom{\indexProsumer}{\indexRiskMeasure} > \valueProsumerStochasticRealization{\indexProsumer}{\indexRiskMeasure}{\indexScenarioDominance}) \;
    \iff \; & \underset{\indexScenario=1}{\overset{\nbScenario}{\sum}} \pi_{\indexScenario} b^\indexProsumer_{\indexScenario,\indexScenarioDominance} \leq a_{\indexScenarioDominance}} & \forall \indexScenarioDominance \in \Omega.
\end{align*}
To linearize the definition of $b^\indexProsumer_{\indexScenario,\indexScenarioDominance}$, we rely on big-M constraint: 
\begin{align*}
    & b^\indexProsumer_{\indexScenario,\indexScenarioDominance} \in \na{0,1} & \forall \indexScenarioDominance \in \Omega, \forall \indexScenario \in \Omega\\
    & \utilityProsumerStochasticRealization{\indexProsumer}{\indexRiskMeasure}{\indexScenario} - \valueProsumerStochasticRealization{\indexProsumer}{\indexRiskMeasure}{\indexScenarioDominance} \leq M b^\indexProsumer_{\indexScenario,\indexScenarioDominance} & \forall \indexScenarioDominance \in \Omega, \forall \indexScenario \in \Omega.
\end{align*}
\end{proof}

\subsection{Increasing convex dominance constraint model}
\label{ann:ic_dominance_constraints}

The increasing convex dominance constraints \eqref{eq:acceptability_seconddc}, is based on \cite{CGS2009}.
\begin{lem}
    In problem $\modelAggregatorStochastic{ }{\indexRiskMeasure}$, the acceptability constraint  $\utilityProsumerStochasticRandom{\indexProsumer}{\indexRiskMeasure}\preceq_{(ic)} \; \valueProsumerStochasticRandom{\indexProsumer}{\indexRiskMeasure}$ can be modeled with: 
    \begin{subequations}
    \begin{align}
        &s^{\indexProsumer}_{\indexScenario,\indexScenarioDominance} \geq 0  & \forall \indexScenarioDominance \in [ \nbScenario], \; \forall \indexScenario \in \Omega
        \\
        & s^{\indexProsumer}_{\indexScenario,\indexScenarioDominance} \geq \utilityProsumerStochasticRealization{\indexProsumer}{\indexRiskMeasure}{\indexScenario} - \valueProsumerStochasticRealization{\indexProsumer}{\indexRiskMeasure}{\indexScenarioDominance}  & \forall \indexScenarioDominance \in [ \nbScenario], \; \forall \indexScenario \in \Omega \\
        & \underset{\indexScenario=1}{\overset{\nbScenario}{\sum}} \pi_{\indexScenario} s^{\indexProsumer}_{\indexScenario,\indexScenarioDominance} \leq a^{ic}_{\indexScenarioDominance} & \forall \indexScenarioDominance \in [ \nbScenario].
    \end{align}
    \end{subequations}
    We denote $a^{ic}_{\indexScenarioDominance} := \EE \nc{(\valueProsumerStochasticRandom{\indexProsumer}{\indexRiskMeasure} -\valueProsumerStochasticRealization{\indexProsumer}{\indexRiskMeasure}{\indexScenarioDominance})^+}$.
\end{lem}

\begin{proof}
    
As in \cref{ann:1_dominance_constraints}, we know that $\valueProsumerStochasticRandom{\indexProsumer}{\indexRiskMeasure}$ follows a discrete distribution with realizations $\valueProsumerStochasticRealization{\indexProsumer}{\indexRiskMeasure}{\indexScenarioDominance})$ for $\indexScenarioDominance \in \Omega$. 
Then,
\begin{align*}
    \utilityProsumerStochasticRandom{\indexProsumer}{\indexRiskMeasure}\preceq_{(ic)} \; \valueProsumerStochasticRandom{\indexProsumer}{\indexRiskMeasure}
    \iff \quad & \EE\bc{\, ( \utilityProsumerStochasticRandom{\indexProsumer}{\indexRiskMeasure} - \eta)^+\, }\, \leq \EE\bc{\, ( \valueProsumerStochasticRandom{\indexProsumer}{\indexRiskMeasure} - \eta)^+ \,}
     & \forall \eta \in \RR \\
    \iff \quad & \EE\bc{\, ( \utilityProsumerStochasticRandom{\indexProsumer}{\indexRiskMeasure} - \valueProsumerStochasticRealization{\indexProsumer}{\indexRiskMeasure}{\indexScenarioDominance})^+\, }\, \leq \EE\bc{\, ( \valueProsumerStochasticRandom{\indexProsumer}{\indexRiskMeasure} - \valueProsumerStochasticRealization{\indexProsumer}{\indexRiskMeasure}{\indexScenarioDominance})^+ \,} & \forall \indexScenarioDominance \in \Omega. 
\end{align*}

We introduce positive variables $s^{\indexProsumer}_{\indexScenario,\indexScenarioDominance} = (\utilityProsumerStochasticRealization{\indexProsumer}{\indexRiskMeasure}{\indexScenario} - \valueProsumerStochasticRealization{\indexProsumer}{\indexRiskMeasure}{\indexScenarioDominance})^+$, for $\eta \in \Omega$.
Thus, we can model the increasing convex dominance constraints as: follows
\begin{align*}
    & \bgp{\EE\bc{\, ( \utilityProsumerStochasticRandom{\indexProsumer}{\indexRiskMeasure} - \valueProsumerStochasticRealization{\indexProsumer}{\indexRiskMeasure}{\indexScenarioDominance})^+\, }\, \leq \EE\bc{\, ( \valueProsumerStochasticRandom{\indexProsumer}{\indexRiskMeasure} - \valueProsumerStochasticRealization{\indexProsumer}{\indexRiskMeasure}{\indexScenarioDominance})^+ \,} \; \iff \; \sum_{\indexScenario=1}^{\nbScenario} \pi_{\indexScenario} s^{\indexProsumer}_{\indexScenario,\indexScenarioDominance} \leq a^{ic}_{\indexScenarioDominance}} & \forall \indexScenarioDominance \in [ \nbScenario].
\end{align*}

\end{proof}

\end{appendix}

\end{document}